\documentclass[11pt]{amsart}
\usepackage{amssymb,amscd}
\usepackage[cmtip,matrix,arrow]{xy}

\newtheorem{theorem}{Theorem}[section]
\newtheorem{corollary}[theorem]{Corollary}
\newtheorem{proposition}[theorem]{Proposition}

\newtheorem{lemma}[theorem]{Lemma}

\theoremstyle{remark}
\newtheorem{remark}[theorem]{Remark}
\newtheorem{definition}[theorem]{Definition}
\newtheorem{remarks}[theorem]{Remarks}
\newtheorem{example}[theorem]{Example}
\newtheorem{examples}[theorem]{Examples}
\newcommand\A{\mathcal{A}}
\newcommand\B{\mathcal{B}}
\newcommand\be{\begin{equation}\label}
\newcommand\ee{\end{equation}}
\newcommand\M{\mathcal{M}}

\newcommand{\N}{\mathcal{N}}
\newcommand{\R}{\mathbb{R}}
\newcommand{\C}{\mathbb{C}}
\newcommand{\F}{\mathbb{F}}

\newcommand{\Z}{\mathbb{Z}}

\renewcommand{\P}{\mathcal{P}}

\newcommand\lie[1]{\mathfrak{#1}}
\renewcommand{\k}{\lie{k}}
\newcommand{\h}{\lie{h}}

\newcommand{\g}{\lie{g}}

\newcommand{\on}{\operatorname}

\newcommand{\Ad}{ \on{Ad} } 
\newcommand{\ad}{ \on{ad} }

\newcommand{\End}{ \on{End} } 

\newcommand{\Hom}{ \on{Hom}} 
\renewcommand{\ker}{ \on{ker}}

\newcommand\dirac{/\kern-1.2ex\partial} 
\newcommand\qu{/\kern-.7ex/} 

\newcommand{\lra}{\longrightarrow}
\newcommand{\hra}{\hookrightarrow}

\renewcommand{\d}{{\on{d}}}
\newcommand{\ol}{\overline}

\newcommand\sig{\sigma}
\newcommand\eps{\epsilon}
\newcommand\Om{\Omega}

\newcommand{\f}{\frac}

\newcommand{\hh}{{\textstyle \f{1}{2}}}
\newcommand{\ti}{\tilde}

\newcommand\beqn{\begin{equation}}      
\newcommand\eeqn{\end{equation}}      
\newcommand{\ca}{\mathcal}

\newcommand{\mf}{\mathfrak}
\newcommand{\beq}{\begin{eqnarray*}}
\newcommand{\eeq}{\end{eqnarray*}}

\begin{document}

\title[]{Equivariant cohomology and\\ the 
Maurer-Cartan equation}

\author{A. Alekseev}
\address{University of Geneva, Section of Mathematics,
2-4 rue du Li\`evre, 1211 Gen\`eve 24, Switzerland}
\email{alekseev@math.unige.ch}

\author{E. Meinrenken}
\address{University of Toronto, Department of Mathematics,
100 St George Street, Toronto, Ontario M5S3G3, Canada }
\email{mein@math.toronto.edu}

\date{\today}

\begin{abstract}
{Let $G$ be a compact, connected Lie group, acting smoothly on a
manifold $M$. In their 1998 paper, Goresky-Kottwitz-MacPherson described a
small Cartan model for the equivariant cohomology of $M$,
quasi-isomorphic to the standard (large) Cartan complex of equivariant
differential forms. In this paper, we construct an explicit cochain
map from the small Cartan model into the large Cartan model,
intertwining the $(S\g^*)_{\on{inv}}$-module structures and inducing
an isomorphism in cohomology. The construction involves the solution 
of a remarkable inhomogeneous Maurer-Cartan equation. This solution 
has further applications to the theory of transgression in the Weil 
algebra, and to the Chevalley-Koszul theory of the cohomology of 
principal bundles.  

\vskip 0.3cm


}
\end{abstract}

\subjclass{ }
\maketitle
\tableofcontents

\section{Introduction}\label{sec:intro}
Let $G$ be a compact connected Lie group of rank $l$, 
and $\pi\colon P\to B$ a principal
$G$-bundle with connection. Choose a connection on $P$.
By the Chern-Weil construction, the de Rham complex 
$\Om(B)$ of differential forms on the base 
becomes a module for the  algebra $(S\g^*)_{\on{inv}}$ of invariant 
polynomials. Consider the corresponding Koszul complex, 
\begin{equation}\label{eq:AAA}
\Om(B) \otimes(\wedge\g^*)_{\on{inv}} ,\ \  \d\otimes 1+
\sum_j p^j\otimes\iota(c_j) ,
\end{equation}
where $c_1,\ldots,c_l$ are the primitive generators of
$(\wedge\g)_{\on{inv}}$, and $p^1,\ldots,p^l$ are the
generators of $(S\g^*)_{\on{inv}}$, corresponding to the dual basis
$c^j\in (\wedge\g^*)_{\on{inv}}$ by Chevalley's transgression theorem.
It is a classical result of Chevalley and Koszul
\cite{ko:tr1,gr:co3}, that the complex \eqref{eq:AAA} is
quasi-isomorphic to the complex $\Om(P)_{\on{inv}}$ of 
invariant forms on the total space. 

Goresky-Kottwitz-MacPherson in \cite{gor:eq} described a similar 
``small'' model for the equivariant de Rham cohomology of any 
$G$-manifold $M$. Recall that the standard Cartan model for the equivariant 
de Rham cohomology of $M$ is the complex
\begin{equation}\label{eq:cart1}
 (S\g^*\otimes\Om(M))_{\on{inv}},\,\,1\otimes \d-\sum_a v^a \otimes\iota(e_a),
\end{equation}
where $e_a\in \g$ is a basis, and $v^a\in S\g^*$ are the generators of the 
symmetric algebra given by the dual basis. By contrast, the small Cartan model 
introduced in  \cite{gor:eq} involves only {\em invariant} differential 
forms: 
\begin{equation}\label{eq:cart2}
 (S\g^*)_{\on{inv}}\otimes\Om(M)_{\on{inv}},\ \ 
1\otimes \d-\sum_j p^j \otimes\iota(c_j).
\end{equation}
The goal of the present paper is the construction of an
explicit cochain map from the small Cartan model \eqref{eq:cart2} into
the Cartan model \eqref{eq:cart1}, commuting with the
$(S\g^*)_{\on{inv}}$-module structure and inducing an isomorphism in
cohomology. In more detail, we will construct a nilpotent even
element 
$$f\in (S\g^*\otimes\wedge\g)_{\on{inv}},$$ 
such that the natural inclusion from \eqref{eq:cart2} to
\eqref{eq:cart1}, followed by a `twist' by the operator $e^{\iota(f)}$
(letting $\wedge\g$ act by contraction), is the desired cochain
map. As it turns out, these properties are equivalent to the following 
Maurer-Cartan equation for the element $f$, 
\begin{equation} \label{eq:!!}
\partial f+\hh [f,f]_{\wedge\g}=\sum_j p^j \otimes c_j-\sum_a v^a \otimes e_a. 
\end{equation}
Here $\partial$ is the Lie algebra boundary operator on $\wedge\g$, 
and $[\cdot,\cdot]_{\wedge\g}$ the Schouten bracket, both extended to 
the algebra $S\g^*\otimes\wedge\g$ in the natural way.
The main discovery of this paper is that this equation does, in
fact, have a solution.

The original argument in \cite{gor:eq} for the equivalence of the two
Cartan models was based on the Koszul duality between differential
$(S\g^*)_{\on{inv}}$-modules and $(\wedge\g)_{\on{inv}}$-modules.  In
this approach, one has to show that \eqref{eq:AAA} is quasi-isomorphic
to $\Om(P)_{\on{inv}}$ not just as a differential space, but also as a
differential $(\wedge\g)_{\on{inv}}$-module. As pointed out in
Mazszyk-Weber \cite{mas:ko}, the proof of this fact in \cite{gor:eq}
contains an error (the proposed map is not a cochain map) 
. 
However, the alternative argument in \cite{mas:ko}
is incorrect as well (the cochain map given there does not respect the
$(\wedge\g)_{\on{inv}}$-module structure). As we will explain in this 
paper, the desired quasi-isomorphism of differential 
$(\wedge\g)_{\on{inv}}$-modules may be constructed by once again 
employing the twist $e^{\iota(f)}$. 

The element $f$ also provides a new point of view on transgression.
Let $W\g=S\g^* \otimes \wedge \g^*$ be the Weil algebra equipped with
the Weil differential $\d^W$. Then, $e^{\iota(f)}$ is an operator
acting on $W\g$. We show that for any primitive element $c^j \in (\wedge \g^*)_{\on{inv}}$
one has
$$
\d^W \left( e^{\iota(f)} (1\otimes c^j) \right) = p^j \otimes 1.
$$
That is, $e^{\iota(f)} (1\otimes c^j)$ is a cochain of transgression
for the polynomial $p^j$.  

We would like to point out some recent references related to this
work. In \cite{all:go}, Allday-Puppe prove a version of a conjecture
of Goresky-Kottwitz-MacPherson that the small Cartan model may be
replaced with an even smaller `Hirsch-Brown' model,
$(S\g^*)_{\on{inv}}\otimes H(M)$, where the differential is
constructed from cohomology operations over $M$.  In a different
direction, Franz \cite{fra:ko1,fra:ko2} has introduced small
Cartan-type models for the equivariant cohomology with {\em integer}
coefficients. Huebschmann \cite{hue:ho,hue:re} obtained
small Cartan models  using homological perturbation techniques. 

The organization of the paper is as follows. In Section
\ref{sec:prel}, we collect some formulas for the Lie algebra boundary
and coboundary operators, and review Chevalley's theory of
transgression in the Weil algebra. In Section \ref{sec:weil} we
consider the problem of finding a cochain map from the Koszul algebra
over the space of primitive elements into the Weil algebra.  This
naturally leads to the above Maurer-Cartan equation. We prove that the
Maurer-Cartan equation admits a solution, which is unique up to
`gauge transformations'. In the subsequent Sections \ref{sec:cartan}
and \ref{sec:chevalley} we apply our results to the small Cartan and
the Chevalley-Koszul complexes. Finally, in the appendix we explain
how the two complexes may be viewed as special cases of a
more general complex due to Halperin \cite{gr:co3}, and how they 
are related by Koszul duality \cite{gor:eq}.

Throughout, we will work in the algebraic context of differential
spaces, over any field $\F$ of characteristic zero.  The applications
to manifolds are obtained as special cases for $\F=\R$, working with
complexes of differential forms.
\vskip.1in

\noindent{\bf Acknowledgements.} We are grateful to M. Franz for
explaining his work and for very useful suggestions. We would like to
thank C. Allday for valuable discussions.  Research of A.A. was
supported in part by the Swiss National Science Foundation. Research
of E.M.  was supported in part by the Natural Sciences and Engineering
Research Council of Canada.

\section{Preliminaries}\label{sec:prel}
In this Section we recall some basic results (due largely to
Chevalley, Hopf, and Koszul) concerning the structure of the invariant
subspace of the symmetric and exterior algebra over any reductive Lie
algebra.  For more details, see \cite{ko:ho,ko:cl,gr:co3}. 

\subsection{Graded vector spaces}
Throughout this paper, $\F$ will denote a field of characteristic $0$.
We will frequently encounter graded vector spaces
$V=\bigoplus_{i\in\Z}V^i$ over $\F$.  Such vector spaces form a
category $\on{GR}_\F$, with morphisms the linear maps preserving
degree. Given a graded vector space $V$, we denote by $V[k]$ the
vector space $V$ with the shifted grading $V[k]^i=V^{k+i}$. A linear
map $V\to W$ between graded vector spaces has degree $k$ if it defines
a morphism $V\to W[k]$.  The tensor product $V\otimes W$ of two graded
vector spaces carries a grading 
$(V\otimes W)^i=\bigoplus_{r+s=i} V^r\otimes W^s$. 
Define the commutativity isomorphism $V\otimes W\to
W\otimes V$ by $v\otimes w\mapsto (-1)^{|v||w|}w\otimes v$, where
$|\cdot|$ denote the degree of a homogeneous element.  Together with the 
obvious associativity isomorphism, $U\otimes (V\otimes W)\to (U\otimes
V)\otimes W$, this makes $\on{GR}_\F$ into a tensor category. One can therefore
consider its algebra objects (called graded algebras), Lie algebra objects
(called graded Lie algebras), and so forth. The commutativity isomorphism
encodes the \emph{super sign convention}: For instance, if
$\A=\bigoplus_{i\in\Z}\A^i$ is a graded algebra, we denote by
$[\cdot,\cdot]$ the (graded) commutator
$$ [x,y]=xy-(-1)^{|x||y|}yx.$$
(This makes $\A$ into a graded Lie algebra; $\A$ is called commutative if the bracket is trivial.) 
A derivation of $\A$  is a linear map $\partial\in \End(\A)$ such that 
\[ [\partial,\eps(x)]=\eps(\partial x)\]
where $\eps\colon \A\to \on{End}(\A)$ is given by left multiplication. 
(Later, we will usually omit $\eps$ from the notation.)
Similarly, one defines derivations of graded Lie algebras. 

\subsection{Lie algebra homology and cohomology}
Let $\g$ be a Lie algebra over $\F$. In order to avoid confusion with
commutators, the Lie bracket will be denoted $[\cdot,\cdot]_\g$. For
any $\g$-module $\M$,  the operator
corresponding to $\xi\in\g$ is denoted $L^\M(\xi)$, or simply 
$L(\xi)$ if $\M$ is clear from the context.

Consider the exterior powers of the adjoint and coadjoint 
representations, with gradings 
\[(\wedge\g)^{-i}=\wedge^i\g,\ \ (\wedge\g^*)^i=\wedge^i\g^*.\]
For $\xi\in\g$ we denote by 
\[ \eps(\xi)\in\End(\wedge\g), \ \iota(\xi)\in\End(\wedge\g^*)\]
the operators of exterior multiplication and contraction.  Note that both 
of these operators have degree $-1$, hence they extend to homomorphisms of graded algebras, 
\begin{equation}
\eps\colon\wedge\g\to \End(\wedge\g),\ \ 
\iota\colon\wedge\g\to \End(\wedge\g^*)\label{eq:iota}.\end{equation}
Dually, for $\mu\in\g^*$ we define operators of degree $+1$,
$\iota^*(\mu)\in \End(\wedge\g),\ \eps^*(\mu)\in \End(\wedge\g^*)$.
Recall Koszul's formulas for the Lie algebra differentials
$\d\in\End(\wedge\g^*)$ and $\partial\in\End(\wedge\g)$, 
\begin{equation}\label{eq:koszuld}
\d=\hh \sum_a \eps^*(e^a) L(e_a),\ \  
\partial=-\hh \sum_a L(e_a) \iota^*(e^a)
\end{equation}
where $e_a\in\g$ and $e^a\in\g^*$ are dual bases. 
Both of these operators square to $0$ and are $\g$-equivariant. 
The differential $\d$ is a derivation of $\wedge\g^*$, while 
$\partial$ is a coderivation of the natural coproduct 
on $\wedge\g$. On the other hand, the interaction of $\partial$ 
with the product on $\wedge\g$ 
is given by the formula (cf. \cite[p.178]{gr:co3})
\begin{equation}\label{eq:product2}
\partial(f\wedge g)=\partial f\wedge g+(-1)^{|f|}f\wedge\partial g
+(-1)^{|f|}[f,g]_{\wedge\g}
\end{equation}
where $[\cdot,\cdot]_{\wedge\g}$ is the {\em Schouten bracket}, 
\begin{equation}\label{eq:schouten}
\begin{split}
[f,g]_{\wedge\g}
&=-\sum_a L(e_a) f\wedge \iota^*(e^a) g.
\end{split}
\end{equation}
The Schouten bracket makes $(\wedge\g)[1]$ into a graded Lie algebra, 
with $\g$ as a Lie subalgebra. The differential $\partial$ 
is a derivation of the bracket, so that  
$(\wedge\g)[1]$ is a differential graded Lie algebra. 
The center of $(\wedge\g)[1]$ is the invariant subspace
$(\wedge\g)_{\on{inv}}[1]$, with the zero differential.  

We will need the following generalization of Cartan's formula 
$[\d,\iota(\xi)]=L(\xi)$ for $\xi\in\g$:
\begin{lemma}\label{lem:comm3}
For any $f\in\wedge\g$, 
\begin{equation}\label{eq:comm3}
 [\d,\iota(f)]=-\iota(\partial f)+\sum_a \iota(\iota^*(e^a)f)\,
 L(e_a).\end{equation}
\end{lemma}
\begin{proof}
The proof is by induction on the degree of $f$, the case 
$|f|=1$ being Cartan's formula.  Suppose 
$f=\xi\wedge g$ where
$|g|=|f|-1$ and $\xi\in\g$. By induction, we may assume that the formula holds
for $g$. Thus  
\[ \begin{split}
[\d,\iota(\xi\wedge g)]&=L(\xi)\iota(g)-\iota(\xi) [\d,\iota(g)]\\
&=\iota(L(\xi)g)+\iota(g)L(\xi)+\iota(\xi)\iota(\partial g)
-\sum_a \iota(\xi) \iota(\iota^*(e^a)g)L(e_a).
\end{split}
\]
The second term can be written as $\iota(\iota^*(e^a)\xi\wedge g)L(e_a)$, 
which combines with the fourth term to $\sum_a\iota(\iota^*(e^a)f)L(e_a)$.
The first and third term add to $-\iota(\partial f)$ since 
$\partial(\xi\wedge g)=-\xi\wedge\partial g -L(\xi)g$. 
\end{proof}

Let $(\wedge\g)^-=\bigoplus_{i>0}\wedge^i\g$. Since any 
$f\in(\wedge\g)^-_{\on{even}}$ is nilpotent, the exponential 
$e^f=\sum_{n=0}^\infty\f{1}{n!}f^n$ is given by a finite series.
\begin{lemma}\label{lem:twist1}
For any $f\in(\wedge\g)^-_{\on{even}}$, 
\begin{equation}
\label{eq:comm2}
 e^{-\iota(f)}\circ \d\circ e^{\iota(f)}
=\d-\iota(\partial f+\hh [f,f]_{\wedge\g})+\sum_a \iota(\iota^*(e^a)f)L(e_a).
\end{equation}
\end{lemma}
\begin{proof}
We may write $e^{-\iota(f)}\circ \d\circ e^{\iota(f)}$ as a sum, 
\[ 
\Ad(e^{-\iota(f)})\d=\d+[\d,\iota(f)]+
\hh \big[[\d,\iota(f)],\iota(f)\big]+\cdots.
\]
The first commutator $[\d,\iota(f)]$ was computed in \eqref{eq:comm3}. 
The next commutator is
\[  \big[[\d,\iota(f)],\iota(f)\big]=
\sum_a \iota\big(\iota^*(e^a)f\wedge L(e_a)f\big)=
-\iota([f,f]_{\wedge\g}), \]
and all higher commutators vanish. 
\end{proof}

\begin{lemma}\label{lem:nice}
For any $f\in(\wedge\g)^-_{\on{even}}$, 
\begin{equation}\label{eq:nice}
 e^{-f}\partial(e^f)=\partial f+\hh [f,f]_{\wedge\g}.
\end{equation}
\end{lemma}
\begin{proof}
Applying \eqref{eq:comm3} to the element $e^f$, we find 
\[\begin{split}
\d\circ \iota(e^f)&=
\iota(e^f)\circ \d-\iota(\partial e^f)
+\sum_a \iota(\iota^*(e^a)e^f)L(e_a)\\
&=\iota(e^f)\circ \big(\d-\iota(e^{-f}\partial e^f)
+\sum_a \iota(\iota^*(e^a)f)L(e_a)\big).\end{split}\]
Equation \eqref{eq:nice} follows by comparing this formula with 
\eqref{eq:comm2}.
\end{proof}

Below we will use these formulas in slightly greater generality: 
Suppose $\B$ is a commutative, evenly graded algebra. 
We use the same notation for the contraction operation \eqref{eq:iota} ,  
the differentials 
\eqref{eq:koszuld}, and the Schouten bracket \eqref{eq:schouten}, 
and for their extensions to $\B\otimes \wedge\g$, $\B\otimes \wedge\g^*$, 
by $\B$-linearity. Then Equations \eqref{eq:comm2} and \eqref{eq:nice} 
hold for any even element $f\in \B\otimes (\wedge\g)^-$, by the same proof. 

\subsection{Hodge theory on $\wedge\g$}\label{subsec:hodge}
Suppose now that $\g$ is a reductive Lie algebra. Then the 
projection onto invariants $\wedge\g\to (\wedge\g)_{\on{inv}}$ 
is a homotopy equivalence for the differential $\partial$, with homotopy 
inverse the inclusion \cite[p.189]{gr:co3}. Recall 
the construction of a homotopy operator, using Hodge theory
with respect to an invariant scalar product $B$ on $\g$. Let 
$B^\flat\colon\g\to \g^*$ and $B^\sharp\colon\g^*\to \g$ 
denote the isomorphisms defined by $B$, and let
$$ \on{Cas}_\g=\sum_a B^\sharp(e^a)e_a\in U(\g)_{\on{inv}}$$
be the quadratic Casimir operator. Let $\on{Cas}_\g^\wedge$ denote the 
operator on $\wedge\g$ defined by $\on{Cas}_\g$ via the adjoint 
representation. Then $\wedge\g=\on{im}\on{Cas}^\wedge_\g
\oplus \on{ker}\on{Cas}^\wedge_\g$, and $\on{ker}\on{Cas}^\wedge_\g$
is the invariant subspace $(\wedge\g)_{\on{inv}}$.
The isomorphisms 
$B^\flat$ and $B^\sharp$ extend to the 
exterior algebras, and in particular the Lie algebra differential $\d$ 
on $\wedge\g^*$ defines a differential on $\wedge\g$
$$\delta=B^\sharp\circ \d\circ B^\flat\in \on{Der}(\wedge\g).$$
The corresponding Hodge Laplacian 
$$ \ca{L}=\delta\partial +\partial\delta\in \End(\wedge\g)$$
has degree $0$, and equals 
$-\f{1}{2}\on{Cas}_\g^\wedge$  \cite[Equation (94)]{ko:cl}.
Hodge theory shows that 
$\on{im}\ca{L}=\on{im}\delta \oplus \on{im}\partial$, 
and one obtains the direct sum decomposition \cite[Proposition 22]{ko:cl},
\begin{equation}\label{eq:hodgedec}
 \wedge\g=\on{im}\delta \oplus \on{im}\partial \oplus (\wedge\g)_{\on{inv}}.
\end{equation}
Let $\ca{G}\in \End(\wedge\g)$ denote the Green's
operator, i.e. $\ker\ca{G}=\ker{\ca{L}}$ and
$\ca{L}\ca{G}=\ca{G}\ca{L}=I-\Pi$ where $\Pi\colon \wedge\g\to
(\wedge\g)_{\on{inv}}$ is the projection defined by the
splitting. Then $\ca{G}$ has degree $0$, and 
$\ca{S}=\ca{G}\delta=\delta \ca{G} $ is the desired
homotopy operator:
$$
[\ca{S},\partial]=[\ca{G}\delta,\partial]=\ca{G}[\delta,\partial]=\ca{G}\ca{L}=I-\Pi.$$
Using that $\g$ is reductive, the differential $\partial$ 
may be written $\partial=-\hh \sum_a  \iota^*(e^a)\ L(e_a)$, and in 
particular vanishes on invariants. Hence one obtains an isomorphism of 
vector spaces, $H(\wedge\g,\partial)=(\wedge\g)_{\on{inv}}$. Similarly, 
the inclusion $(\wedge\g^*)_{\on{inv}}\hra \wedge\g^*$ 
defines an isomorphism of algebras 
$H(\wedge\g^*,\d)=(\wedge\g^*)_{\on{inv}}$. 

\subsection{Primitive elements}
Since $\g$ is reductive, the pairing between $\wedge\g$ and
$\wedge\g^*$ restricts to a non-degenerate pairing between $(\wedge
\g)_{\on{inv}}$ and $(\wedge \g^*)_{\on{inv}}$. As a consequence, the
algebra structure on $(\wedge \g^*)_{\on{inv}}$ induces a coalgebra
structure on $(\wedge \g)_{\on{inv}}$, and the algebra structure 
on $(\wedge \g)_{\on{inv}}$ induces a coproduct on $(\wedge \g^*)_{\on{inv}}$.
The coproduct
and product are compatible in both cases, turning $(\wedge
\g)_{\on{inv}}$ and $(\wedge \g^*)_{\on{inv}}$ into commutative graded
Hopf algebras.  Recall that an element $x$ of a graded coalgebra is
called {\em primitive} if it has the property,
$$ \Delta(x)=x\otimes 1+1\otimes x,$$
where $\Delta$ is the coproduct. Let $\P,\P^*$ denote the graded
subspaces of primitive elements in $(\wedge \g)_{\on{inv}},(\wedge
\g^*)_{\on{inv}}$, respectively. It can be shown \cite[page
206]{gr:co3} that the pairing between $(\wedge \g)_{\on{inv}},(\wedge
\g^*)_{\on{inv}}$ restricts to a non-degenerate pairing between $\P$
and $\P^*$, so that indeed $\P^*$ is the dual space to $\P$. 
By results of Hopf and Samelson, the elements
in $\P,\P^*$ all have odd degree, and the inclusion maps extend to
graded Hopf algebra isomorphisms
$$ \wedge \P\stackrel{\cong}{\lra} (\wedge \g)_{\on{inv}},\ \ \wedge \P^*\stackrel{\cong}{\lra}
(\wedge \g^*)_{\on{inv}}.$$
This means in particular that the operator of contraction by any $c\in \P$
is a {\em derivation} of $(\wedge\g^*)_{\on{inv}}$, even though it is
not of course a derivation of $\wedge\g^*$.

\subsection{Transgression}
Let $S\g^*$ be the symmetric algebra over $\g^*$, with grading 
$$(S\g^*)^{2i}=S^i\g^*,\ \ (S\g^*)^{2i+1}=0.$$ 
Let $\ti{\P}=\P[-1]$ be the 
evenly graded vector
space, obtained by lowering the grading of $\P$ by $1$, 
and dually $\ti{\P}^*=\P^*[1]$. 
There is a canonical isomorphism of graded algebras, due to Koszul 
and Chevalley (see \cite[page 242]{gr:co3})
\[  S\ti{\P}^*\cong (S\g^*)_{\on{inv}}. \]
Let us review Chevalley's construction \cite[page 363]{gr:co3}
of this isomorphism, using transgression in the Weil algebra
$$
 W\g=S\g^*\otimes \wedge\g^*.
$$
Fix a basis $e_a$ of $\g$, with dual basis $e^a$, and let
$y^a\in\wedge^1\g^*$ and $v^a\in S^1\g^*$ be the corresponding
generators of the exterior and symmetric algebra. The Weil
differential $\d^W$ is the derivation of $W\g$ given by the formula,
\footnote{From now on, we identify the elements of any algebra with
the corresponding operator of left multiplication on the algebra.}
\begin{equation}\label{eq:weildiff}
 \d^W=\sum_a y^a L^W(e_a) -\d^\wedge +\sum_a v^a \iota(e_a).
\end{equation}
Here $L^W(e_a)=L^S(e_a)+L^\wedge(e_a)$ are the generators for the
$\g$-action, while $\d^\wedge,\iota$ are the Lie algebra
differential and the contraction operator, acting on the second
factor, $\wedge\g^*$. The Weil algebra is acyclic, and so is the
invariant subalgebra $(W\g)_{\on{inv}}$. The subspace 
$(S\g^*)_{\on{inv}}\subset W\g$ consists of cocycles for the 
Weil differential. Hence, by acyclicity, any element in 
$(S^+\g^*)_{\on{inv}}$ is exact.

An odd element $x\in(W\g)_{\on{inv}}$ is called a {\em cochain of
transgression} if its differential $\d^Wx$ lies in
$(S\g^*)_{\on{inv}}$.  It is called a {\em distinguished cochain of
transgression} if, in addition, $\iota(c)x\in \F$ for all $c\in P$.
The space of cochains of transgression is denoted $\ca{T}$, and its
subspace of distinguished cochains of transgression
$\ca{T}_{\on{dist}}$. The image of $\ca{T}\subset W\g$ 
under the projection $W\g\to \wedge\g^*$ (defined by 
the augmentation map $S\g^*\to \F$) is exactly 
the subspace $\P^*\subset (\wedge\g^*)_{\on{inv}}$ 
of primitive elements. Together with the map
$\ca{T}\to (S\g^*)_{\on{inv}}, \,x\mapsto \d^W x$ this fits into a
commutative diagram,
$$ \xymatrix{ &\ca{T}\ar[ld]\ar[d]&\\
(S^+\g^*)_{\on{inv}}\ar[r]_{\ \ \ \ \tau} &\P^*&
\!\!\!\!\!\!\!\!\!\!\!\!\!\!\subset (\wedge\g^*)_{\on{inv}}}
$$
The map $\tau\colon(S^+\g^*)_{\on{inv}}\to \P^*$ has degree $-1$ and is 
referred to as {\em transgression}. Its kernel is the annihilator of 
the ideal $\big((S^+\g)_{\on{inv}}\big)^2\subset (S^+\g)_{\on{inv}}$ 
of decomposable elements.
As it turns out \cite[page 239]{gr:co3}, the map from 
$\ca{T}_{\on{dist}}$ to $(\wedge\g^*)_{\on{inv}}$ is still onto $\P^*$, 
and there is a {\em unique} map $\gamma\colon\P^*\to (S\g^*)_{\on{inv}}$ of degree $1$ 
such that the following diagram commutes: 
$$ \xymatrix{ &\ \ \ca{T}_{\on{dist}}\ar[ld]\ar[d]&\\
(S^+\g^*)_{\on{inv}}&\P^*\ar[l]^{\ \ \ \ \ \  \gamma}&
\!\!\!\!\!\!\!\!\!\!\!\!\!\!\!\!\!\subset (\wedge\g^*)_{\on{inv}}}
$$
The map $\gamma$ identifies $\ti{\P}^*$ as a subspace of
$(S\g^*)_{\on{inv}}$.  By Chevalley's theorem, the inclusion map
extends to an algebra isomorphism, $S\ti{\P}^*\cong
(S\g^*)_{\on{inv}}$.

Below we will also need another description of the 
space $\P^*$ of primitive elements. Let 
\begin{equation}\label{eq:sigma}
\varsigma\colon S\g^*\to
\wedge\g^*\end{equation}
be the homomorphism of graded algebras, given on
$S^1\g^*=\g^*$ by the Lie algebra differential. View 
$\wedge\g^*$ as a module for the subalgebra $\on{im}(\varsigma)$, 
and let $\on{im}(\varsigma)\g^*$ be the submodule generated by $\g^*$. 
\begin{lemma}\label{lem:kost}
The space of primitive elements in $(\wedge\g^*)_{\on{inv}}$ is 
the invariant subspace of $\on{im}(\varsigma)\g^*$:
$$\P^*=(\on{im}(\varsigma)\g^*)_{\on{inv}}.$$
\end{lemma}
\begin{proof}
It is well-known (cf. \cite[page 233]{gr:co3} or \cite[Equation (261)]{ko:cl})
that for any polynomial $p\in (S^i\g^*)_{\on{inv}}$ the
element $\tau(p)$ is given, up to a multiplicative constant, by
$$ \sum_a \,\varsigma(\iota^S(e_a)p)\wedge e^a \in (\wedge^{2i-1}\g^*)_{\on{inv}}.$$
Here  $\iota^S(\xi)$ is the derivation of $S\g^*$, given on $S^1\g^*=\g^*$ 
by the natural pairing. (Put differently, $\iota^S(\xi)p$ is the 
derivative of the polynomial $p$ in the direction of $\xi$.)  This gives the 
inclusion $\P^* \subset (\on{im}(\varsigma)\g^*)_{\on{inv}}$. 
Since $\g$ is reductive, the map
$$\Hom_\g(\g,\on{im}(\varsigma))\cong
(\on{im}(\varsigma)\otimes\g^*)_{\on{inv}} \to
(\on{im}(\varsigma)\g^*)_{\on{inv}}$$
given by wedge product is onto. According to Kostant
\cite[Equation (263)]{ko:cl}, the multiplicity of the adjoint
representation in $\on{im}(\varsigma)$ equals $\on{rank}(\g)=\dim \P^*$. We
conclude
$$\dim (\on{im}(\varsigma)\g^*)_{\on{inv}}
\le \dim \Hom_\g(\g,\on{im}(\varsigma))
=\dim \P^*.$$ 
\end{proof}
 
\section{A canonical cochain of transgression}\label{sec:weil}
\subsection{The inclusion $K(\P)\hra W\g$}
The action of $\wedge\g$ on $\wedge\g^*$ by contractions $\iota$
extends to an action (still denoted $\iota$) of the graded algebra
$S\g^*\otimes \wedge\g$ on the Weil algebra $W\g=S\g^*\otimes
\wedge\g^*$. Lemma \ref{lem:comm3} shows that on the invariant
subspace $(W\g)_{\on{inv}}$, the action of $(S\g^*)_{\on{inv}}\otimes
(\wedge\g)_{\on{inv}}$ commutes with the differential. That is,
$(W\g)_{\on{inv}}$ is a differential graded module over the algebra 
$(S\g^*)_{\on{inv}}\otimes (\wedge\g)_{\on{inv}}$.

Let $c^j\in \P^*$ and $c_j\in \P$ denote dual (homogeneous) bases for
the primitive subspaces.  Let $p^j=\gamma(c^j)\in \ti{\P}^*$ denote the
corresponding generators of $(S\g^*)_{\on{inv}}$.  The Koszul algebra
of $\P$ is the tensor product
$$K(\P)=S\ti{\P}^*\otimes\wedge \P^*,$$
with differential $\d^K=\sum_j p^j\otimes \iota(c_j)$. Thus $\d^K$
vanishes on $S\ti{\P}^*$, and takes the generators $c^j$ of $\wedge
\P^*$ to the corresponding generators $p^j$ of $S\ti{\P}^*$.  
The Koszul algebra is a differential graded modules for 
the algebra
$$S\ti{\P}^*\otimes\wedge \P,$$
where the first factor acts by multiplication and the second factor
acts by contraction.  The natural inclusion 
$K(P)\hra (W\g)_{\on{inv}}$ is compatible with the module structures,
but is not a cochain map. However, we have the following result:
\begin{theorem}\label{th:weil}
There is an injective homomorphism of graded vector spaces 
$$ \Phi\colon K(\P)\to (W\g)_{\on{inv}},$$
with the following properties: 
\begin{enumerate}
\item \label{a}
$\Phi$ is a cochain map,  
\item\label{b}
$\Phi$ is a homomorphism of modules over 
$S\ti{\P}^*\otimes\wedge \P\cong  (S\g^*)_{\on{inv}}\otimes
(\wedge\g)_{\on{inv}}$, 
\item\label{c}
$\Phi$ is unital, that is, it takes the unit of $K(P)$ to 
the unit of $(W\g)_{\on{inv}}$.
\end{enumerate}
\end{theorem}
As a direct consequence of Theorem \ref{th:weil}, 
we have:
\begin{corollary}\label{cor:trans}
The restriction of the map $\Phi$ to $\P^*=\wedge^1\P^*\subset K(P)$
fits into a commutative diagram, 
$$ \xymatrix{ &\ \
\ca{T}_{\on{dist}}\ar[ld]&\\(S^+\g^*)_{\on{inv}}&\P^*\ar[u]_\Phi\ar[l]^{\
\ \ \ \ \ \gamma}& \!\!\!\!\!\!\!\!\!\!\!\!\!\!\!\!\!\subset
(\wedge\g^*)_{\on{inv}}}$$
\end{corollary}
\begin{proof}
The element $\Phi(c^j)\in (W\g)_{\on{inv}}$ satisfies
$$ \d^W \Phi(c^j)=\Phi(\d^K c^j)=\Phi(p^j)=p^j\Phi(1)=p^j=\gamma(c^j).$$
Since furthermore $\iota(c_i)\Phi(c^j)=\Phi(\iota(c_i)c^j)=\delta_i^j$, it 
follows that $\Phi(c^j)$ is a distinguished cochain of transgression.
\end{proof}

\begin{remark}\label{rem:chevhom}
On the other hand, let $\ti{c}^j\in
(W\g)_{\on{inv}}$ be distinguished cochains of transgression 
extending $c^j$, and consider the algebra homomorphism, 
\[ \Phi'\colon K(P)\to (W\g)_{\on{inv}},\ \ \ \Phi'(p^j)=p^j,\ \Phi'(c^j)=\ti{c}^j.\]
Then $\Phi'$ is a homomorphism of differential $S\ti{P}^*$-algebras, but it does not intertwine the $\wedge
P$-actions (unless $\g$ is Abelian). Indeed, recall that the operators $\iota(c_j)$ 
are derivations of $K(P)$ but not of $(W\g)_{\on{inv}}$. 
\end{remark}

\subsection{The map $\Phi$}\label{subsec:twists}
In this Section we determine the most general form of a map
$\Phi$ which satisfies conditions \eqref{b} and \eqref{c}, and 
reduce the condition \eqref{a} to a Maurer-Cartan type equation.
Note that for any even element 
$f\in (S\g^*\otimes(\wedge\g)^-)_{\on{inv}}$, the map 
\begin{equation}\label{eq:anyf}
\Phi\colon  
S\ti{\P}^*\otimes \wedge \P^* \cong 
(S\g^*)_{\on{inv}} \otimes (\wedge \g^*)_{\on{inv}}
\hra (W\g)_{\on{inv}} \stackrel{e^{\iota(f)}}{\lra} (W\g)_{\on{inv}} 
\end{equation}
satisfies properties \eqref{b} and \eqref{c}. Furthermore, 
$\Phi$ preserves degrees if and only if $f$ has degree $0$.
The following converse was pointed out to us by M. Franz: 
\begin{lemma}
Any even linear map $\Phi\colon K(P)\to (W\g)_{\on{inv}}$ satisfying
conditions \eqref{b} and \eqref{c} is of the form \eqref{eq:anyf}, for
a unique even element $f\in (S\g^*\otimes(\wedge\g)^-)_{\on{inv}}$. 
Moreover, 
$\Phi$ preserves degrees if and only if $f$ has degree $0$. 
\end{lemma}
\begin{proof}
Observe that any element of $K(P)$ is obtained from the volume element
$c^1\cdots c^l$ by the action of $S\ti{\P}^*\otimes\wedge \P$.  Hence,
a map $\Phi$ satisfying \eqref{b} is uniquely determined by
$\Phi(c^1\cdots c^l)$. 
Similarly, any element in $W\g$ is obtained
from $c^1\cdots c^l$ (now viewed as a volume element in $\wedge\g^*$)
by the action of $S\g^*\otimes\wedge\g$.  Let $F\in
(S\g^*\otimes\wedge\g)_{\on{inv}}$ be the unique element such that
$$\iota(F)(c^1\cdots c^l)=\Phi(c^1\cdots c^l).$$
Then $\Phi$ is a composition of the inclusion map 
$K(P)\hra (W\g)_{\on{inv}}$ with
$\iota(F)$. The image of the unit $1\in K(P)$ under this map 
equals the component of $F$ in
$(S\g^*\otimes \wedge^0\g)_{\on{inv}}\cong (S\g^*)_{\on{inv}}$. 
Hence, by (c) this component must be equal to $1$.
Equivalently, $F=e^f$ where $f\in (S\g^*\otimes (\wedge\g)^-)_{\on{inv}}$ is 
given by 
$$f=\log(F)=- \sum_n (1-F)^n/n.$$
(This is well-defined since $1-F$ is nilpotent.) 
If $\Phi$ preserves the grading, $F$ as defined above has degree $0$, 
hence also $f=\log(F)$ has degree $0$.
\end{proof}

Our next task is to arrange that $\Phi$ is a cochain map.

\begin{proposition}\label{prop:weildiff}
For any even element $f\in (S\g^*\otimes(\wedge\g)^-)_{\on{inv}}$, the conjugate of 
the Weil differential by $e^{-\iota(f)}$ is given by the formula, 
\begin{equation}\label{eq:weilconj}
 \Ad(e^{-\iota(f)})\d^W=\d^W+\iota\big(\partial f+\hh[f,f]_{\wedge\g}\big)
+\sum_a \iota\big(\iota^*(e^a) f\big) \,L^S(e_a).
\end{equation}
\end{proposition}
\begin{proof}
The first two terms in Formula \eqref{eq:weildiff} 
for the Weil differential contribute 
\[ \begin{split}
 \Ad\big(e^{-\iota(f)}\big)\Big(\sum_a y^a L^W(e_a)\Big)
&=\sum_a y^a L^W(e_a)+\sum_a \iota\big(\iota^*(e^a)f\big)L^W(e_a),\\
\Ad\big(e^{-\iota(f)}\big)(-\d^\wedge)
&= -\d^\wedge+\iota(\partial f+\hh [f,f]_{\wedge\g})-
\sum_a \iota(\iota^*(e^a)f)L^\wedge(e_a)
\end{split}\]
(using $\Ad(e^{-\iota(f)})y^a=y^a+\iota(\iota^*(e^a)f)$ and 
\eqref{eq:comm2}), while the last term 
$\sum_a v^a \iota(e_a)$ in \eqref{eq:weildiff}
commutes with the action of $e^{-\iota(f)}$. 
Equation \eqref{eq:weilconj} follows. 
\end{proof} 
By \eqref{eq:weilconj} and the formula \eqref{eq:weildiff} for 
the Weil differential, we obtain
\[ \Ad(e^{-\iota(f)})\d^W=\sum_a v^a \otimes \iota(e_a)
+\iota\big(\partial f+\hh[f,f]_{\wedge\g}\big)+\cdots,\]
where the dots indicate terms vanishing on 
$(S\g^*)_{\on{inv}}\otimes (\wedge\g^*)_{\on{inv}}$. 
Hence $\Phi=e^{\iota(f)}$
will give the desired cochain map $K(\P)\to (W\g)_{\on{inv}}$, 
provided $f\in (S\g^*\otimes(\wedge\g)^-)_{\on{inv}}$ has degree $0$ 
and solves the inhomogeneous Maurer-Cartan equation, 
\begin{equation}\label{eq:basic}
\partial f+\hh [f,f]_{\wedge\g}=
\sum_j p^j \otimes c_j-\sum_a v^a \otimes e_a.
\end{equation}

\begin{theorem}\label{th:exists}
The inhomogeneous Maurer-Cartan equation \eqref{eq:basic}
has a (canonical) solution $f\in (S\g^*\otimes(\wedge\g)^-)_{\on{inv}}$ 
of degree $0$. In fact, $Z=\sum_j p^j \otimes c_j$ is the
\emph{only} element in $(S\g^*)_{\on{inv}}
\otimes(\wedge\g)_{\on{inv}}$ for which the equation 
\[\partial f+\hh [f,f]_{\wedge\g}+\sum_a v^a \otimes e_a
=Z\]
admits such a solution.
\end{theorem}
Corollary \ref{cor:trans} may now be restated as the assertion that 
for any even, nilpotent solution $f$ of Equation \eqref{eq:basic} and any 
$c^j\in \P^*$, the element
$\ti{c}^j=e^{\iota (f)}c^j$ is a distinguished cochain of
transgression, with $\d^W \ti{c}^j=p^j$.

The proof of Theorem \ref{th:exists} will be given at the end of Section
\ref{subsec:solution}, after some preparations.  As a consequence of
Theorem \ref{th:exists}, we recover Chevalley's correspondence between
primitive generators of $(\wedge\g)_{\on{inv}}$ and generators of
$(S\g^*)_{\on{inv}}$ from a rather unexpected angle: It takes the form of a  
solvability
condition for an inhomogeneous Maurer-Cartan equation in a differential 
graded Lie algebra. 

Let us introduce the following notation, 
\begin{equation}\label{eq:ourcase}
\begin{split}
\k&=\bigoplus_{i\le 0} \k^i,\ \ \  \k^{i}=(S\g^*\otimes \wedge^{1-i}\g)_{\on{inv}},\\
\mf{l}&=\bigoplus_{i\le 0} \mf{l}^i,\ \ \mf{l}^i=(S\g^*)_{\on{inv}}\otimes 
(\wedge^{1-i}\g)_{\on{inv}},\\
X&=-\sum_a v^a\otimes e_a.
\end{split}
\end{equation}
As a first step toward solving \eqref{eq:basic}, we will look for solutions 
of $\partial f+\hh [f,f]_{\k}=X$ \emph{modulo} $\mf{l}$. This will 
be done in Section \ref{subsec:general} below. 
We will need the following fact:
\begin{lemma}
The element $X=-\sum_a v^a \otimes e_a$ is a cocycle, 
contained in the center $\mf{z}$. 
\end{lemma}
\begin{proof}
It is clear that $X$ is a cocycle, since $\partial$ vanishes on 
$\g\subset \wedge\g$. Furthermore, using invariance of elements $\phi\in\k$ 
under the diagonal action, 
\[ [X,\phi]_{\wedge\g}=-\sum_a v^a L^\wedge(e_a)\phi=\sum_a v^a L^S(e_a) \phi=0.\]
Here we used the fact that  the derivation $\sum_a v^a L^S(e_a)$ of 
$S\g^*$ is zero.
\end{proof}
Observe that $\mf{l}$ is contained in the center $\mf{z}$ of $\k$, and 
that $\partial$ vanishes on $\mf{l}$. Furthermore, since 
$(\wedge\g)_{\on{inv}}\hra \wedge\g$ is a $\g$-equivariant 
homotopy equivalence, the inclusion of $\mf{l}$ into $\k$
is a homotopy equivalence.

\subsection{Solutions of a Maurer-Cartan equation}\label{subsec:general}
It is convenient to place \eqref{eq:basic} into a more general
framework. Let $\k=\bigoplus_{i\leq 0} \k^i$ be a differential graded
Lie algebra with a differential $\partial$ of degree $+1$. We assume
that the grading is bounded below, that is $\k^i =0$ for $i<\!< 0$.
Denote by $U(\k)=\bigoplus_{i\leq 0} U(\k)^i$ the enveloping algebra,
with grading induced by the grading of $\k$ (that is, 
${\rm deg} \, (x_1\cdots x_r) = i_1+\ldots +i_r$
for $x_j\in \k^{i_j}$),  
and by $\ol{U}(\k)=\prod_{i\leq 0}
U(\k)^i$ the degree completion of $U(\k)$.  
Elements of $\ol{U}(\k)$ are infinite
series $a=\sum_{i\leq 0} a_i$ with $a_i \in U(\k)^i$. The differential
$\partial$ extends to $\ol{U}(\k)$ as a derivation of the product.
Write $\k^-=\bigoplus_{i<0} \k^i$.  Its even part $\k^-_{\on{even}}$
is an (ordinary) nilpotent Lie algebra.\footnote{In the previous
Section, the labels `even' and `odd' referred to the grading on
$\wedge\g$.  Since \eqref{eq:ourcase} involves the shifted grading
$\wedge\g[1]$, the roles of `even' and `odd' are now reversed.}  There
is a well-defined exponential map 
$\exp: \k^-_{\on{even}} \rightarrow
\ol{U}(\k), s \mapsto \exp(s)=\sum_N\f{1}{N!} s^N$.  It is a 1-1 map,
and its image in $\ol{U}(\k)$ is a group, with product given by the
Campbell-Hausdorff formula.
One has the following well-known formula for the (right) {\em
Maurer-Cartan form}, 
\begin{equation}\label{eq:MCform}
(\partial\exp(s))\exp(-s)=j^R(\ad_s)\partial s,\ \ \ s\in
\k^-_{\on{even}},
\end{equation}
where $j^R(z)=\f{e^z-1}{z}$ and $\ad_s=[s,\cdot]_\k$. 
The Maurer-Cartan form enters into the formula for the {\em gauge action}
of $\exp(\k^-_{\on{even}})$ on $\k^-_{\on{odd}}$,
\begin{equation}\label{eq:gauge}
\exp(s).f=e^{\ad_s}f-j^R(\ad_s)\partial s.
\end{equation}
The {\em curvature} 
$\partial f+\hh [f,f]_\k$ of any $f\in \k^-_{\on{odd}}$ 
transforms under the adjoint 
representation: If $\ti{f}= \exp(s).f$, then 
\begin{equation} \label{eq:curvature}
\partial \ti{f}+\hh [\ti{f},\ti{f}]_\k
=e^{\ad_s}\big(\partial f+\hh [f,f]_\k\big).
\end{equation}
\begin{theorem}\label{th:general}
Let $(\k,\partial)$ be a differential graded Lie algebra, with 
$\k^i=0$ for $i<\!<0$ and $\k^i=0$ for $i>0$. Assume there is a 
subspace $\mf{l}$ of the center $\mf{z}$ of $\k$, 
such that $\partial$ vanishes on $\mf{l}$ and such that 
the inclusion $\mf{l}\hra \k$ induces an isomorphism in cohomology. 
Then:
\begin{enumerate}
\item \label{it:a}
For any even central element  $X\in \mf{z}$ with  $\partial X=0$, 
the set of solutions $f\in \k^-_{\on{odd}}$ of the equation
\begin{equation}\label{eq:basic2}
 \partial f+\hh[f,f]_{\k}=X\ \mod \mf{l}
\end{equation}
is a homogeneous space for the group 
$\exp(\k^-_{\on{even}})\times\mf{l}^-_{\on{odd}}$, 
where the first factor acts by gauge transformations 
and the second factor by translations. 
\item \label{it:c}
The difference
\[ \partial f+\hh[f,f]_{\k}-X\in \mf{l}  \]
is independent of the solution $f$. 
\end{enumerate}
\end{theorem}
\begin{proof}
For any solution $f$ of \eqref{eq:basic2} 
the curvature $\partial f+\hh [f,f]_\k$ is in the center of $\k$. It 
is therefore invariant under gauge transformations of $f$
(see \eqref{eq:curvature}). 
On the other hand, the curvature is also invariant under 
the translation action of $\mf{l}^-_{\on{odd}}$
(since element of $\mf{l}$ are central cocycles by assumption). 
Hence  \eqref{it:c} follows from \eqref{it:a}. 

To prove \eqref{it:a} we have to show that the solution space is
non-empty and that the action of
$\exp(\k^-_{\on{even}})\times\mf{l}^-_{\on{odd}}$ is transitive.  It
suffices to prove these statements for the quotient
$\k/\mf{l}$. Equivalently, we may (and will) assume for the rest of
this proof that $\mf{l}=0$, hence $H(\k)=0$.
Write $f=f_1+f_3+\ldots$ with $f_i\in\k^{-i}$, and similarly 
$X=X_0+X_2+\ldots$ with $X_i\in \k^{-i}$. 
Then \eqref{eq:basic2} is equivalent to a system of equations, 
\begin{equation}
\tag{$A_N$}\ \ \  \ \ \ \ \partial f_N =-\hh \sum_{i+j=N-1}[f_i,f_j]_{\k}+X_{N-1},\  \ 
\end{equation}
where $N=1,3,\ldots$. Equation ($A_1$) reads $\partial f_1 =X_0$.  
It admits a solution since $\partial X_0=0$ and since $H(\k)=0$. 
Suppose by induction that we have found solutions 
$f_i$ for the system of equations $(A_i)$ 
up to $i=N-2$, and consider $(A_N)$. We have, 
$$ \partial\big(-\hh \sum_{i+j=N-1}[f_i,f_j]_{\k}+X_{N-1}\big)=
-\sum_{i+j=N-1}[\partial f_i,f_j]_{\k} =-\sum_{i+j=N-1}[X_{i-1},f_j]_{\k} =0.$$
Here we have used the Jacobi identity for $\k$, and the assumption that 
$X$ is central. Since $H(\k)=0$, it follows that $(A_N)$ admits a solution. 
This proves the existence part of the theorem. To show uniqueness 
up to gauge transformations, suppose $f\in \k^-_{\on{odd}}$ is a solution of 
\eqref{eq:basic2}. Then the derivation 
$\nabla=\partial+\ad_f$ is again a differential on $\k$:
$$ \nabla^2=\ad(\partial f+\hh[f,f]_{\k})=\ad(X)=0.$$
Given $r\in \k^-_{\on{odd}}$, the sum $f+r$ is a 
solution of \eqref{eq:basic2} if and only if 
\begin{equation}\label{eq:r}
 \nabla r+\hh [r,r]_{\k}=0.
\end{equation}
We must show that $f+r$ is gauge equivalent to $f$.
The formula \eqref{eq:gauge} for gauge transformations of $f$ can be 
written
\begin{equation}\label{eq:gauge1}
 \exp(s).f=f-j^R(\ad_s)\nabla s,
\end{equation}
where we used \eqref{eq:MCform} and the identity
$e^{\ad_s}f=f+j^R(\ad_s)[s,f]_\k$. Hence, we have to prove
that 
\begin{equation}\label{eq:general form}
r=-j^R(\ad_s)\nabla s
\end{equation} 
for some $s=s_2+s_4+\cdots$ with $s_i\in \k^{-i}$. 
Write 
\[ s_{\{ N\}}=s_2+s_4+\cdots + s_N,\ \ \ 
r_{\{ N\}}=-j^R(\ad_{s_{\{ N\}}})\nabla s_{\{ N\}}.
\]
Suppose we have found $s_2,s_4,\ldots,s_N$, such that $q_{\{
N\}}=r-r_{\{ N\}}$ is contained in $\k^{-(N+1)}\oplus
\k^{-(N+3)}+\cdots$. Since $r_{\{
N\}}$ solves the Maurer-Cartan equation \eqref{eq:r}, 
$$ \nabla q_{\{ N\}}+\hh
[q_{\{ N\}},q_{\{ N\}}]_{\k}+[r_{\{ N\}},q_{\{ N\}}]_{\k}=0.$$
By construction, the left hand side lies in $\k^{-N}\oplus 
\k^{-(N+2)}+\cdots$. Moreover, the only contribution to the component 
in $\k^{-N}$ comes from $\partial q_{\{ N\}}$. It follows that the component 
of $q_{\{ N\}}$ in $\k^{-(N+1)}$ is closed, hence exact (since $H(\k)=0$). 
Choose $s_{N+2}\in \k^{-(N+2)}$ such that the component of 
$q_{\{ N\}}-\partial s_{N+2}$ in $\k^{-(N+1)}$ is zero. Letting
$s_{\{ N+2\}}:=s_{\{ N\}}+s_{N+2}$, we achieve 
$q_{\{ N+2\}}\in \k^{-(N+3)}\oplus \k^{-(N+5)}+\cdots$.
Hence, by induction we obtain the desired element $s$. 
\end{proof}
Assume that $\ca{S}\colon \k\to \k$ is a 
homotopy operator, i.e. $\ca{S}$ has degree $-1$ and 
$[\ca{S},\partial]=I-\Pi$ where
$\Pi$ is a projection operator onto $\mf{l}$. Then we can write down an explicit 
solution to \eqref{eq:basic2}, by the following recursion formula: 
\begin{equation}\label{eq:recursion}
 f_N=\ca{S} \big(-\hh \sum_{i+j=N-1}[f_i,f_j]_{\k}+X_{N-1}\big),\ \ 
N=1,3,\ldots
\end{equation}
\subsection{Solution of Equation \eqref{eq:basic}}\label{subsec:solution}
We now return to our original problem, Equation \eqref{eq:basic}. 
Choose an invariant scalar product $B$ on $\g$, and let   
$\ca{S}=\delta\ca{G}$ be the homotopy operator defined by Hodge theory 
(Section \ref{subsec:hodge}). 
Let $f=f_1+f_3+\ldots$ be the solution {\em modulo}
$(S\g^*)_{\on{inv}}\otimes (\wedge\g)_{\on{inv}}$,  
given by the recursion formula \eqref{eq:recursion}. 
\begin{lemma}
The element $f\in (S\g^*\otimes (\wedge\g)^-)_{\on{inv}}$ has total 
degree $0$.
\end{lemma}
Here \emph{total degree} refers to our original grading on 
the algebra $S\g^*\otimes \wedge\g$, in contrast to the 
Lie algebra grading used in \eqref{eq:ourcase}.
\begin{proof}
The element $X=\sum_a v^a\otimes e_a$ has total degree $1$, while 
the homotopy operator $\ca{S}$ has total degree $-1$. 
Thus $f_1=\ca{S}(X)$ has total degree $0$. Since the Schouten bracket of any 
two elements of total degree $0$ has total degree $-1$, the recursion formula 
\eqref{eq:recursion} shows that all $f_N$ have total degree $0$. 
\end{proof}

To complete the proof of Theorem \ref{th:exists}, 
it remains to identify the `error term' 
\begin{equation}\label{eq:Z}
Z=\partial f+\hh [f,f]_{\wedge\g}+\sum_a v^a \otimes e_a
\in (S\g^*)_{\on{inv}}\otimes (\wedge\g)_{\on{inv}}.\end{equation}
Recall that by Theorem \ref{th:general}\eqref{it:c}, $Z$ is independent the choice
of solution $f$.

Using the isomorphism $B^\flat\colon \g\to \g^*$, the map
\eqref{eq:sigma} translates into an algebra homomorphism,
\[\zeta\colon S\g\to \wedge\g.\]
extending the map $\delta\colon\g\to \wedge^2\g$.  
View $\wedge\g$ as a module for the subalgebra $\on{im}(\zeta)$, and let
$\on{im}(\zeta)\g$ denote the submodule generated by
$\g=\wedge^1\g$.  Notice that
\[
\begin{split}
\iota^*(B^\flat(\xi))&\colon  \on{im}(\zeta)\to \on{im}(\zeta)\g,\\
         L(\xi)&\colon  \on{im}(\zeta)\to \on{im}(\zeta),\\
         \delta&\colon  \on{im}(\zeta)\g\to  \on{im}(\zeta),\\
       \partial&\colon \on{im}(\zeta)\to \on{im}(\zeta)\g.
\end{split}
\]
It follows that 
$\on{im}(\zeta)$ and $\on{im}(\zeta)\g$ are both 
invariant under $\ca{L}=[\delta,\partial]$, 
and that the Schouten bracket of any two elements 
in $\on{im}(\zeta)$ is contained in $\on{im}(\zeta)\g$. 
\begin{remark}\label{rem:zeta}
The subspace $\on{im}(\zeta)$ does not depend on the choice of
the invariant scalar product $B$. If $\g$ is simple, this follows 
since $B$ is unique up to a multiplicative constant in this case, and 
$\delta$ just scales by that constant. For the general case, it
suffices to observe that if $\g=\bigoplus \g_i$ with scalar product 
$B=\bigoplus B_i$, then the subalgebra $\on{im}(\zeta)$ is 
generated by the images of the differentials $\delta_i$, obtained from 
the Lie algebra differential on $\g_i$ by the isomorphism
$B_i^\flat\colon \g_i\to \g_i^*$.
\end{remark}

\begin{lemma}\label{lem:zeta}
The solution $f\in (S\g^*\otimes \wedge\g)_{\on{inv}}$ defined by
the homotopy operator 
$\ca{S}=\delta \ca{G}$ is contained in $(S\g^*\otimes
\on{im}(\zeta))_{\on{inv}}$.  Hence,
$$ Z \in 
(S\g^*)_{\on{inv}}\otimes (\on{im}(\zeta)\g)_{\on{inv}}
=(S\g^*)_{\on{inv}}\otimes \P.$$
\end{lemma}
\begin{proof}
We use the notation from Section \ref{subsec:hodge}. 
Since $\on{im}(\zeta)$ is invariant under $\ca{L}$, 
it is also invariant under the Green's operator $\ca{G}$.
Hence, $\ca{S}\colon (S\g^*\otimes \on{im}(\zeta)\g)_{\on{inv}}\to  
(S\g^*\otimes \on{im}(\zeta))_{\on{inv}}$. Obviously, 
$X=-\sum_a v^a \otimes e_a \in (S\g^*\otimes \on{im}(\zeta)\g)_{\on{inv}}$.
An induction based on the recursive definition 
\eqref{eq:recursion} therefore shows that 
each $f_n$ is in $(S\g^*\otimes \on{im}(\zeta))_{\on{inv}}$. This proves the first claim. The properties of $\zeta$ show 
$Z\in S\g^*\otimes \on{im}(\zeta)\g$, while on the other hand
$Z\in (S\g^*)_{\on{inv}}\otimes (\wedge\g)_{\on{inv}}$. 
Finally, $(\on{im}(\zeta)\g)_{\on{inv}}=\P$ by Lemma \ref{lem:kost}.
\end{proof}
\begin{proof}[End of Proof of Theorem \ref{th:exists}]
Using the lemma we may write $Z =\sum_j q^j \otimes c_j$, 
where the $c_j\in (\wedge\g)_{\on{inv}}$ are a 
homogeneous basis of $\P$, while 
the $q^j \in (S\g^*)_{\on{inv}}$ are invariant polynomials. 
To determine the elements $q^j$, let us once again consider 
$e^{\iota(f)}$ as an operator on $(W\g)_{\on{inv}}$.  
By Proposition \ref{prop:weildiff}
we have 
$$\Ad(e^{-\iota(f)})(\d^W)=\sum_j q^j \otimes \iota(c_j)+\ldots=\iota(Z)+\ldots
$$
where $\ldots$ 
are terms vanishing on ${(S\g^*)_{\on{inv}}\otimes (\wedge\g^*)_{\on{inv}}}$.
Apply this result to $c^k\in \P^*=(\wedge\g^*)_{\on{inv}}$.
The element $\ti{c}^k=e^{\iota(f)}c^k \in (W\g)_{\on{inv}}$ satisfies, 
\[ \begin{split}
\d^W\,\ti{c}^k
&=\d^W e^{\iota(f)}c^k =e^{\iota(f)}\iota(Z)c^k
=e^{\iota(f)} q^k=q^k,\\
\iota(c_j)\ti{c}^k&=\iota(c_j)e^{\iota(f)}c^k
=e^{\iota(f)} \delta_j^k=\delta_j^k.\end{split}\]
Thus, the $\ti{c}^k$ are distinguished cochains of transgression. In 
particular, $q^k=\d^W\,\ti{c}^k=\gamma(c^k)=p^k$. This concludes the 
proof of Theorem \ref{th:exists}. 
\end{proof}

According to Theorem \ref{th:general}, the solution $f$ of the
Maurer-Cartan equation \eqref{eq:basic} is unique up to gauge 
transformation and translation by elements in
$(S\g^*)_{\on{inv}}\otimes (\wedge\g)_{\on{inv}}$. The special 
solution found above is singled out by the `gauge fixing'. 
More precisely: 
\begin{proposition}
The Maurer-Cartan equation \eqref{eq:basic} admits a unique 
$\delta$-exact solution $f\in (S\g^*\otimes(\wedge\g)^-)_{\on{inv}}$,
given explicitly by the recursion formula \eqref{eq:recursion} with $\ca{S}=\delta\ca{G}$. 
This solution does not depend on the choice of invariant 
scalar product $B$ (even though $\delta$ does). 
\end{proposition}
\begin{proof}
  Recall that the Maurer-Cartan equation \eqref{eq:basic} is
  equivalent to a system of equations of the form ($A_N$) for
  $f=f_1+f_3+\cdots\in (S\g^*\otimes(\wedge\g)^-)_{\on{inv}}$, 
\[ 
 \partial f_N =-\hh \sum_{i+j=N-1}[f_i,f_j]_{\k}+\Big( \sum_j p^j
  \otimes c_j-\sum_a v^a\otimes e_a \Big)_{N-1}.
\]
At each stage, the recursion \eqref{eq:recursion} picks out the unique (by the Hodge
  decomposition \eqref{eq:hodgedec}) $\delta$-exact solution.  As
  shown in Lemma \ref{lem:zeta}, this solution $f$ is in fact
  contained in the subspace
  $(S\g^*\otimes\on{im}(\zeta))_{\on{inv}}\subset
  \on{im}(\delta)_{\on{inv}}$.  By Remark \ref{rem:zeta} this subspace
  does not depend on $B$.  Hence $f$ does not depend on $B$.
\end{proof}

The first term $f_1$ in our recursion formula for $f$ can be computed
quite easily. Suppose that $\g$ is simple. 
On $\g\subset \wedge\g$, the operator $\on{Cas}_\g^\wedge$ acts as a scalar,
$(\dim\g)^{-1}\on{tr}_\g(\on{Cas}_\g)$. Hence
$$f_1= -\delta \ca{G}(\sum_a v^a e_a)=\f{2 \dim\g}{\on{tr}_\g(\on{Cas}_\g)}
\sum_a v^a \delta e_a. $$ 
\begin{example}\label{ex:su2} 
Let $\g$ be the three-dimensional Lie algebra
with basis $e_1,e_2,e_3$  and bracket relations, $[e_1,e_2]_\g=e_3$, 
$[e_2,e_3]_\g=e_1$, $[e_3,e_1]_\g=e_2$. Using the inner product $B$ on 
$\g$ for which the $e_a$ are an orthonormal basis, we find 
$\on{tr}_\g(\on{Cas}_\g)=-6$. 
Hence, 
$$f_1=-\delta (\sum_a v^a e_a)
= v^1 \otimes (e_2\wedge e_3)+v^2\otimes (e_3\wedge e_1)+
v^3\otimes (e_1\wedge e_2).$$ 
In this example, higher corrections do not appear, so that  
$f=f_1$. Indeed, taking $c=e_1\wedge e_2\wedge e_3$
as a generator of $(\wedge\g)_{\on{inv}}$, a short calculation shows
$$ \hh [f,f]_{\wedge\g}= p \otimes c,\ \ 
\partial f=-\sum_a v^a e_a,$$
where $p=\sum_a v^a v^a\in (S\g^*)_{\on{inv}}$.  
\end{example}
In all of our applications, the solution $f$ of the Maurer-Cartan
equation enters via its exponential in the algebra
$(S\g^*\otimes\wedge\g)_{\on{inv}}$.  We will now give an explicit
formula for the exponential.
\begin{proposition}\label{eq:special}
The exponential of the solution $f$ is given by the formula, 
$$ e^f=(I+\ca{G}\circ \sum_a v^a \delta e_a)^{-1}(1).$$
\end{proposition}
The inverse is well-defined, since $\ca{G}\circ \sum_a v^a \delta e_a$ 
is nilpotent.
\begin{proof} 
  Let $F=e^f$, and denote by $F_{[k]},\,k=0,2,\ldots$ the component in
  $(S\g^*\otimes\wedge^k \g)_{\on{inv}}$. Since $f$ is $\delta$-exact,
  $F$ is $\delta$-closed, and each $F_{[k]}$ with $k\ge 2$ is
  $\delta$-exact.  Using Lemma \ref{lem:nice},
  Equation \eqref{eq:basic} is equivalent to
  \begin{equation}\label{eq:alt}F_{[0]}=1,\ \  \partial F= Y F,
  \end{equation}
  where $Y=-\sum_a v^a \otimes e_a+\sum_j p^j \otimes c_j$. Applying
  $\delta$, we obtain $\ca{L}F=\delta\partial F=(\delta Y) F$. That
  is, 
  \[ F_{[0]}=1,\ \ \ \ca{L}F_{[k+2]}=(\delta Y) F_{[k]} \] 
  because the Hodge Laplacian $\ca{L}$ preserves the
  exterior algebra degree, while $\delta Y=-\sum_a v^a\delta e_a\in
  (S\g^*\otimes\wedge^2 \g)_{\on{inv}}$. Since $F_{[k+2]}$ is
  $\delta$-exact, $F_{[k+2]}=\ca{G}\ca{L}F_{[k+2]}$. 
  We arrive at the recursion formula
  \[ F_{[0]}=1,\ \ F_{[k+2]}=\ca{G}(\delta Y F_{[k]}),\]
  with solution $F_{[k]}=(\ca{G}\circ \delta Y)^k(1)$. Summing
  $F=\sum_k F_{[k]}$ as a geometric series, the proof is complete.
\end{proof}

\section{The small Cartan complex}\label{sec:cartan}

\subsection{$\g$-differential spaces}\label{subsec:gds}
A \emph{$\g$-differential space} is a differential graded vector 
space $(\M,\d)$, together with linear maps 
\begin{equation}\label{eq:lieder}
 L\colon \g\to \End(\M),\ \ 
\iota\colon \g\to \End(\M)\end{equation}
such that the \emph{Lie derivatives} $L(\xi)$ have 
degree $0$ and the \emph{contractions} $\iota(\xi)$ have degree $-1$, 
and such that the following relations hold:
\begin{equation}
\begin{split} 
[\d,\iota(\xi)]&=L(\xi),\\
[L(\xi),\iota(\xi')]&=\iota([\xi,\xi']_\g),\\ 
[\iota(\xi),\iota(\xi')]&=0.
\end{split}
\end{equation}
These relations and the Jacobi identity imply that the Lie derivatives 
$L(\xi)$ define a representation of $\g$ on $\M$, commuting 
with the differential. 
The motivating example of a $\g$-differential space is the space
$\M=\Om(M)$ of differential forms on a manifold $M$ with 
an action of a Lie group $G$, with \eqref{eq:lieder} the 
Lie derivatives and contractions by the infinitesimal generators
of the action. Another example is $\M=\wedge\g^*$, with $\d$
the Lie algebra differential, $\iota(\xi)$ the usual contraction
operators, and $L(\xi)$ the coadjoint representation.  The Weil algebra
$W\g=S\g^*\otimes\wedge\g^*$ is a $\g$-differential space, with $\d$
the Weil differential $\d^W$, contractions $\iota(\xi)=1\otimes
\iota(\xi)$, and Lie derivatives $L(\xi)=L^W(\xi)$.

For any $\g$-differential space $\M$, one defines the horizontal subspace
$\M_{\on{hor}}=\bigcap_{\xi\in\g} \ker(\iota(\xi))$, the invariant 
subspace $\M_{\on{inv}}=\bigcap_{\xi\in\g} \ker(L(\xi))$, and the 
basic subspace $\M_{\on{basic}}=\M_{\on{hor}}\cap \M_{\on{inv}}$. 
Both $\M_{\on{inv}}$ and $\M_{\on{basic}}$ are stable under $\d$. 

The contraction operators on a $\g$-differential space 
$\M$ extend to an algebra homomorphism 
$\iota\colon\wedge\g\to \End(\M)$. The formulas \eqref{eq:comm3} 
and \eqref{eq:comm2} hold in this greater generality, since their 
proof only relied on the commutation relations between the operators $L(\xi),\iota(\xi),\d$. 
Formula \eqref{eq:comm3} shows that the action of 
$(\wedge\g)_{\on{inv}}$ on $\M_{\on{inv}}$ by contractions 
commutes with the differential. That is, $\M_{\on{inv}}$ is a {\em
differential $(\wedge\g)_{\on{inv}}$-module}.

\subsection{The small Cartan model}
The {\em equivariant cohomology} $H_\g(\M)$ of any $\g$-differential space $\M$
is defined as the cohomology of the Cartan complex
\begin{equation}\label{eq:cartan}
C_\g(\M)=(S\g^*\otimes \M)_{\on{inv}},\ \d_\g=1\otimes \d
-\sum_a v^a\otimes \iota(e_a).
\end{equation}
The algebra $(S\g^*)_{\on{inv}}$ of invariant polynomials acts on 
$C_\g(\M)$ by multiplication, and this action commutes with the differential. 
That is, $C_\g(\M)$ is a {\em differential $(S\g^*)_{\on{inv}}$-module}. 
\begin{definition}
The {\em small Cartan model} for $\M$ is the 
differential $(S\g^*)_{\on{inv}}$-module 
\begin{equation}\label{eq:smallcartan}
 \ti{C}_\g(\M)=(S\g^*)_{\on{inv}}\otimes \M_{\on{inv}},\ \ 
\ti{\d}_\g=1\otimes\d-\sum_j p^j\otimes \iota(c_j).
\end{equation}
\end{definition}
Note that if the Lie algebra $\g$ is Abelian, the Cartan model and the
small Cartan model coincide. In general, we have:
\begin{theorem}\label{th:cartan}
Let $\M$ be any $\g$-differential space. For any solution 
$f\in (S\g^*\otimes(\wedge\g)^-)_{\on{inv}}$ of the Maurer-Cartan equation \eqref{eq:basic}, the composition
\begin{equation}\label{eq:cartanmap}
\ti{C}_\g(\M)\hra C_\g(\M)\stackrel{e^{\iota(f)}}{\lra} C_\g(\M)
\end{equation}
is a homotopy equivalence of differential
$(S\g^*)_{\on{inv}}$-modules.  In particular, it induces an
isomorphism in cohomology, $\ti{H}_\g(\M)\to H_\g(\M)$. 
\end{theorem}
\begin{proof}
By Equations \eqref{eq:comm2} and  \eqref{eq:basic}, the 
operator $e^{\iota(f)}$ on $C_\g(\M)$ takes the Cartan differential
$\d_\g$ into 
\begin{equation}\label{eq:twisteddiff}
\begin{split}
\d_\g':&=e^{-\iota(f)}\circ
\Big(1\otimes \d-\sum_a v^a \otimes \iota(e_a)\Big)\circ e^{\iota(f)}
\\&=1\otimes\d -p^j \otimes \iota(c_j)+\sum_a \iota(\iota^*(e^a)f)L^\M(e_a).\end{split}
\end{equation}
On $(S\g^*)_{\on{inv}}\otimes \M_{\on{inv}}$, the last term on the
right hand side vanishes. This proves that $e^{\iota(f)}$ gives a cochain
map $\ti{C}_\g(\M)\to C_\g(\M)$. 

To construct a homotopy inverse, let $C=C_\g(\M)$ (with differential
$\d_\g'$) and $\ti{C}=\ti{C}_\g(\M)$ (with differential
$\ti{\d}_\g$). We have to find a morphism of differential 
$(S\g^*)_{\on{inv}}$-modules $C\to \ti{C}$ that is homotopy 
inverse to the inclusion, by an $(S\g^*)_{\on{inv}}$-equivariant 
homotopy. Pick an invariant inner product
$B$ on $\g$, and let $\on{Cas}_\g^S$ be the corresponding Casimir operator
on $S\g^*$ (i.e. the image of $\on{Cas}_\g\in U(\g)$ 
under the coadjoint representation). Then $S\g^*$ splits as a direct sum of 
the kernel and image of $\on{Cas}_\g^S$, and the kernel is the space 
of invariants. Let $\ca{L}_0$ denote the restriction of $\on{Cas}_\g\otimes 1$ 
to invariants in $C=(S\g^*\otimes \M)_{\on{inv}}$. Then 
\[ C=\on{ker}(\ca{L}_0)\oplus \on{im}(\ca{L}_0)\]
where $\on{ker}(\ca{L}_0)=\ti{C}$. Let $\Pi_0$ be the projection 
from $C$ onto $\ker(\ca{L}_0)$ along $\on{im}(\ca{L}_0)$, and  
$\ca{G}_0$ the Green's operator, i.e.  $\ca{G}_0\Pi_0=0$ and 
$\ca{L}_0\ca{G}_0=\ca{G}_0\ca{L}_0=1-\Pi_0$. 

Introduce a filtration $0=\M_{(0)}\subset \M_{(1)}\subset \cdots
\subset\M_{(\dim\g+1)}=\M$, where $\M_{(j)}$ is the subspace for which all contractions with
elements in $\wedge^j \g$ are zero.  Since the filtration of $\M$
is $\g$-invariant, it gives rise to a filtration $0=C_{(0)}\subset
C_{(1)}\subset \cdots\subset C_{(\dim\g+1)}=C$. All terms in $\d_\g'$, 
except for the first term $1\otimes\d$, lower the filtration degree by 
at least $1$. Consider the following operator on $C$, 
$$h=-\sum_a L^S(e_a)\otimes \iota(B^\sharp({e}^a))$$
and put $\ca{L}=[\d_\g',h]$. On $C$, 
$$[1\otimes \d,h]=-\sum_a L^S(e_a)\otimes L^\M(B^\sharp({e}^a))
=\on{Cas}_\g\otimes 1=\ca{L}_0.$$
Hence $ \ca{L}=\ca{L}_0+R$ 
where $R$ lowers the filtration degree. From $h\Pi_0=0$ we deduce
$\ca{L}\Pi_0=0$, and therefore
$$ \ca{L}=\ca{L}(1-\Pi_0)=
(1+R\ca{G}_0)\ca{L}_0.$$
But $1+R\ca{G}_0$ is invertible since $R \ca{G}_0$ lowers the 
filtration degree. It follows that 
$C=\on{ker}(\ca{L})\oplus \on{im}(\ca{L})$ 
with  $\on{ker}(\ca{L})=\on{ker}(\ca{L}_0)=\ti{C}$.  
Let $\ca{G}$ denote the Green's operator for the cochain map $\ca{L}$, thus 
$\ca{G}\ca{L}=\ca{L}\ca{G}=1-\Pi$ where $\Pi$ is the projection from 
$C$ onto $\ti{C}$ along $\on{im}(\ca{L})$. Then $H:=h\ca{G}$ satisfies $[\d_\g',H]=1-\Pi$. 
\end{proof}

\begin{remarks}\label{rem:remark}
\begin{enumerate}
\item
The fact that the map \eqref{eq:cartanmap} is a
quasi-isomorphism also follows very quickly from the spectral
sequences for the filtrations
$$ F^j=\bigoplus_{i\ge j} (S^i\g^*\otimes \M)_{\on{inv}},\ \ \ 
\ti{F}^j=\bigoplus_{i\ge j} (S^i\g^*)_{\on{inv}}\otimes \M_{\on{inv}},
$$
provided the complex $\M$ is bounded below (to ensure convergence of
the spectral sequence). Since \eqref{eq:cartanmap} is filtration preserving, it induces a
morphism of spectral sequences, $\phi_r\colon\ti{E}_r\to E_r$. The $0$th
stage for both spectral sequences are the Cartan complexes themselves,
with differential $1\otimes \d$. Hence
$$ \phi_1\colon\ti{E}_1=(S\g^*)_{\on{inv}}\otimes H(\M)\to 
E_1=(S\g^*)_{\on{inv}}\otimes H(\M).$$
Here we have used that 
$$ H(\M_{\on{inv}})=H(\M)_{\on{inv}}=H(\M)$$
since
$L(\xi)=[\iota(\xi),\d]$. The map $\phi_1$ is 
just the identity map, since $e^{\iota(f)}-1$ raises the filtration 
degree by at least two. Since $\phi_1$ is an isomorphism, 
the map in cohomology $\ti{H}_\g(\M)\to H_\g(\M)$ is an isomorphism 
as well.
\item
The $\Z$-grading on $\M$ does not play any role in the proof -- the
results hold more generally for $\g$-differential spaces that are only
$\Z_2$-graded. 
\item
The differential in the small Cartan model is not a derivation for the
obvious product structure. However, it is possible to introduce a new
(non-associative) product, such that the differential is a derivation
and such that the induced product in cohomology is the standard one.
This will be explored in Section \ref{sec:liehom}. It was pointed out 
to us by M. Franz that by the homotopy equivalence between the two 
Cartan models, the algebra structure on the large Cartan model 
gives rise to an $A_\infty$-structure on the small Cartan model. 
The relevant 
machinery is developed in a paper by Gugenheim and Lambe \cite{gug:per}.
\end{enumerate}
\end{remarks}

The small Cartan model is useful provided one has good control over the 
contraction operators $\iota(c_j)$. For instance, it may happen that 
these operators all vanish: 
\begin{example}
Let $M=G/G^\sigma$ be a symmetric space, defined by an {\em inner}
involutive automorphism $\sigma$ of a connected reductive Lie group $G$. 
One example is the
Grassmannian $\on{Gr}_\C(k,n)$ of $k$-planes in $\C^n$, with
$G=U(n)$, and $\sig$ the involution given as conjugation by a diagonal
matrix with entries $(1,\ldots,1,-1,\ldots,-1)$ down the diagonal,
with $k$ plus signs and $n-k$ minus signs.

Recall that $G$-invariant forms on a symmetric space $M$ are
automatically closed, and hence that the cohomology ring of $M$ is
canonically isomorphic to the ring of invariant forms on $M$.  Since $\sig$ is by
assumption an inner automorphism, i.e. $\sig=\Ad_x$ for some $x\in G$,
the induced action on $(\wedge\g)_{\on{inv}}$ is
trivial. In particular, each $c_j$ is $\sig$-invariant.  Write
$\g=\g^\sig\oplus\mf{p}$ where $\mf{p}$ is the $-1$ eigenspace of
$\sig$. Then $\mf{p}$ is identified with the tangent space to $M$ at
the identity coset, and the multi-vector field defined by $c_j$ is
given as the projection of $c_j$ onto $\wedge\mf{p}$. But since $c_j$
has odd degree, and $\sig$ acts as $-1$ on $\wedge^{\on{odd}}\mf{p}$,
it follows that the vector field is just $0$.  That is, $\iota(c_j)=0$
for all $j$. We conclude that any invariant differential form on $M$
is an {\em equivariant} cocycle for the small Cartan model. Applying
the operator $e^{\iota(f)}$, we directly get the equivariant extension
for the standard (large) Cartan model.
\end{example}

\subsection{Dependence on $f$}  \label{sec:choicef}
Returning to the setting of Theorem \ref{th:cartan}, it is natural to
ask to what extent the isomorphism $\ti{H}_\g(\M)\to H_\g(\M)$ depends on the
choice of $f$. Recall that any two solutions $f_0,f_1$ of \eqref{eq:basic} 
are gauge equivalent, up to addition of an even element in 
$\mf{l}=(S\g^*)_{\on{inv}}\otimes (\wedge\g)^-_{\on{inv}}$. 
The action of $(\wedge\g)_{\on{inv}}^-$ on the small Cartan model 
is homotopic to the trivial action: Indeed, let $h_i$ denote the 
derivation of $(S\g^*)_{\on{inv}}$ given by $h_i(p^j)=\delta_i^j$. Then 
$$[\ti{\d}_\g,h_j\otimes 1]=-\iota(c_j),$$
showing that $\iota(c_j)$ is homotopic to $0$. It is therefore
sufficient to consider the case that $f_0,f_1$ are gauge equivalent:
\[ f_1=\exp(s_1).f_0
=e^{\ad_{s_1}}f_0-j^R(\ad_{s_1})\partial s_1\]
for an odd element 
$s_1\in (S\g^*\otimes (\wedge\g)^-)_{\on{inv}}$. 
(Here $\exp$ is the Lie algebra exponential, 
defined as in Section \ref{subsec:general}.) Let $s(t)=ts_1$ 
for $t\in \F$. 
Then $f_0,f_1$ belong to a family of solutions $f(t)$
of \eqref{eq:basic}, depending \emph{polynomially} on $t$ 
(since $s_1$ is nilpotent):
\[f(t)=\exp(s(t)).f_0=e^{\ad_{s(t)}}f_0-
j^R(\ad_{s(t)})\partial s(t).\]
We have, 
\begin{equation}\label{eq:derivative}
 \f{d f}{d t}+\partial s+[f,s]_{\wedge\g}=0.
\end{equation}
Now let $\M$ be any $\g$-differential space. 
\begin{proposition}\label{prop:homotopy}
Let  $\Phi(t)=e^{\iota(f(t))}\colon \ti{C}_\g(\M)\to C_\g(\M)$ be 
the family of cochain maps defined by $f(t)$.  
Then the family of operators
$$H(t)=\Phi(t)\circ \iota(s(t)) \,\colon (S\g^*)_{\on{inv}}\otimes \M_{\on{inv}}\to (S\g^*\otimes
\M)_{\on{inv}},\ 
$$
satisfies, $\f{d \Phi}{d t}=H\circ \ti{\d}_\g+\d_\g \circ H$. Thus,
all $\Phi(t)$ are homotopic as homomorphisms of
$(S\g^*)_{\on{inv}}$-differential modules. 
\end{proposition}
Note that $\Phi(t)$ and $H(t)$ depend polynomially on $t$. 
\begin{proof}
We thank the referee for the following simplified version of 
our original argument. By \eqref{eq:twisteddiff}, we have
\[\d_\g\circ \Phi=\Phi\circ 
\big(\ti{\d}_\g+\sum_a \iota(\iota^*(e^a)f)L^\M(e_a)\big).\]
Furthermore, the commutator of $\sum_a \iota(\iota^*(e^a)f)L^\M(e_a)$ with 
$\iota(s)$ gives $-\iota([f,s]_{\wedge\g})$. 
Hence, using Lemma \ref{lem:comm3} and Equation \eqref{eq:derivative}, 
\[\begin{split}
H\circ \ti{\d}_\g+\d_\g \circ H&=\Phi \circ\iota(s) 
\circ\ti{\d}_\g +\d_\g \circ\Phi \circ\iota(s)\\
&=\Phi \circ \big([\ti{\d}_\g,\iota(s)]
+\sum_a \iota(\iota^*(e^a)f)\circ L^\M(e_a)\circ \iota(s)\big)\\
&=\Phi \circ \big([\ti{\d}_\g,\iota(s)]-\iota([f,s]_{\wedge\g})\big)+\ldots\\
&=-\Phi \circ \iota(\partial s+ [f,s]_{\wedge\g})+\ldots \\
&=\Phi\circ \iota\big(\f{d f}{d t})+\ldots \\
&=\f{d \Phi}{d t}+\ldots,
\end{split}\]
where the dots indicate terms vanishing on $\ti{C}$.
\end{proof}

To summarize, we have shown:
\begin{theorem}
The homomorphism of differential $(S\g^*)_{\on{inv}}$-modules $\Phi\colon 
\ti{C}_\g(\M)\to C_\g(\M)$ is independent of the solution $f$ 
of \eqref{eq:basic},
up to $(S\g^*)_{\on{inv}}$-equivariant homotopy. In particular,
the induced map in cohomology does not depend on the choice of $f$.
\end{theorem}

\section{The Chevalley-Koszul complex}\label{sec:chevalley}
\subsection{$\g$-differential $W\g$-modules}
A \emph{$\g$-differential algebra} is a graded associative algebra $\A$,
together with the structure of a $\g$-differential space in such a way
that the operators $\d,L(\xi),\iota(\xi)$ are derivations for
the product structure. A \emph{$\g$-differential module} for $\A$ is a
$\g$-differential space $\N$, with an $\A$-module structure such that
the action map $\A \otimes \N\to \N$ is a homomorphism of
$\g$-differential spaces. Of particular importance is the case where
$\A=W\g$ is the Weil algebra. If $\N$ is $\g$-differential
$W\g$-module, then the basic subcomplex $\N_{\on{basic}}$ is naturally
a differential $(W\g)_{\on{basic}}=(S\g^*)_{\on{inv}}$-module.
$\g$-differential $W\g$-modules were studied by Guillemin-Sternberg in
\cite{gu:su}, under the name '$W^*$-modules'.
\begin{examples}
\begin{enumerate}
\item
Suppose $\N$ is a commutative $\g$-differential algebra, equipped with an 
algebraic connection $\theta\colon \g^*\to \N$ in the sense of Cartan 
\cite{ca:no}. Recall that the algebraic Chern-Weil homomorphism 
$CW_\theta\colon  W\g\to \N$
is the unique homomorphism of $\g$-differential algebras extending the map  
$\theta$ on $\wedge^1\g^*\subset W\g$. Clearly, $CW_\theta$ gives 
$\N$ the structure of a $\g$-differential $W\g$-module. 
\item
For any $\g$-differential space $\M$, the tensor product
$\N=W\g\otimes \M$ is a $\g$-differential $W\g$-module in the
obvious way. Note that since $\F\hra W\g$ is a homotopy equivalence,
the inclusion $\M\hra \N$ is a quasi-isomorphism of
$\g$-differential spaces.
\end{enumerate}
\end{examples}
As before, we let $y^a$ and $v^a$ denote the generators of $W\g$ 
for a given basis of $\g$.
For any $\g$-differential $W\g$-module $\N$, 
there is a canonical {\em horizontal projection} operator
$$P_{\on{hor}}= 
\prod_a \iota^\N(e_a) y^a\colon \N\to \N_{\on{hor}}.
$$
The operator $P_{\on{hor}}$ is a $\g$-equivariant morphism of 
modules over $S\g^*$. 

We will need the following fact:
\begin{theorem}[Cartan \cite{ca:la}, see also \cite{gu:su,ni:ca}]
\label{th:cargu}
Let $\N$ be a $\g$-differential $W\g$-module. The map 
$ S\g^*\otimes \N\to \N_{\on{hor}},\ p\otimes x\mapsto P_{\on{hor}}(p.x)$
restricts to a \emph{cochain map}, 
\begin{equation}
\label{eq:cart}
C_\g(\N)=(S\g^*\otimes \N)_{\on{inv}}\to \N_{\on{basic}}.
\end{equation}
This cochain map is an $(S\g^*)_{\on{inv}}$-equivariant projection, and 
is homotopy inverse (in the category of
$(S\g^*)_{\on{inv}}$-modules) to the natural inclusion
$\N_{\on{basic}}\to (S\g^*\otimes \N)_{\on{inv}}$.
\end{theorem}
It is useful to note that any $\g$-differential
$W\g$-module $\N$ may be written as a tensor product of $\wedge\g^*$
with an $S\g^*$-module.  
Let $\d_{\on{hor}}=P_{\on{hor}}\circ\d$ is the {\em covariant
derivative} on $\N$, given on horizontal elements $x\in \N_{\on{hor}}$ by
$\d_{\on{hor}}x=(\d-\sum_a y^a\,L^\N(e_a))x.$ 
\begin{proposition}\cite[Section 8.7]{gr:co3}
The action of $\wedge\g^*\subset W\g$
on $\N$ restricts to a $\g$-equivariant isomorphism,
\begin{equation}\label{eq:canonicalisomorphism}
\N_{\on{hor}}\otimes \wedge\g^* \to \N,\ x\otimes\eta\mapsto (-1)^{|\eta||x|}
\eta.x.
\end{equation}
This isomorphism intertwines the action of $W\g=S\g^*\otimes \wedge\g^*$
on $\N_{\on{hor}}\otimes \wedge\g^*$ with the given action on $\N$.
It induces the following differential 
on $\N_{\on{hor}}\otimes\wedge\g^* $, 
\begin{equation}\label{eq:thedifferential}
-1\otimes \d^\wedge+\d_{\on{hor}}\otimes 1
+\sum_a v^a\otimes \iota^\wedge(e_a)+\sum_a (1\otimes y^a)L(e_a).
\end{equation}
Here $L(\xi)=L^\N(\xi)\otimes 1+1\otimes L^\wedge(\xi)$ are the Lie
derivatives for the diagonal action.
\end{proposition}
Note that for $\N=W\g$, one recovers the formula \eqref{eq:weildiff}
for the Weil differential. 

\begin{proof}
Define two $\g$-equivariant degree $0$ endomorphisms of $\N\otimes
\wedge\g^*$,
\[ \alpha=\sum_a \iota^\N(e_a)\otimes y^a,\ \  
\beta=\sum_a y^a\otimes \iota^\wedge(e_a).\]
Since $\alpha,\beta$ are nilpotent, their exponentials are well-defined 
automorphisms of degree $0$. An elementary calculation shows, 
\[ 
\Ad(e^\alpha)(1\otimes\iota^\wedge(\xi))=\iota^\N(\xi)\otimes 1+1\otimes\iota^\wedge(\xi)\\
=\Ad(e^{-\beta})(\iota^\N(\xi)\otimes 1)
\]
Thinking of $\N$ as the joint kernel of the operators $1\otimes\iota^\wedge(\xi)$ 
on $\N\otimes \wedge\g^*$, this gives isomorphisms, 
\begin{equation}\label{eq:isso}  
\N_{\on{hor}}\otimes \wedge\g^*\stackrel {e^{\beta}}{\lra}
(\N\otimes\wedge\g^*)_{\on{hor}} \stackrel{e^{-\alpha}}{\lra}
\N.\end{equation}
On $(\N\otimes\wedge\g^*)_{\on{hor}}$, the operator 
$e^{-\alpha}=\prod_a (1-\iota^\N(e_a)\otimes y^a)$ 
coincides with 
\[ 1\otimes \prod_a (1-y^a\iota^\wedge(e_a))=
1\otimes \prod_a \iota^\wedge(e_a)y^a=1\otimes P^\wedge_{\on{hor}},\] 
where $P_{\on{hor}}\colon\wedge\g^*\to \F$ is the horizontal 
projection map for $\wedge\g^*$. But a moment's 
reflection shows that 
$(1\otimes P^\wedge_{\on{hor}})\circ e^\beta\colon \N_{\on{hor}}\otimes 
\wedge\g^*\to \N$ is exactly the map, 
$x\otimes\eta\mapsto (-1)^{|\eta||x|}\eta.x$. Clearly, this map 
is compatible with the $W\g$-module structures. 

Now let $\d'$ be the differential on $\N_{\on{hor}}\otimes\wedge\g^*$,  
induced via \eqref{eq:isso} from the differential on $\N$, and let 
$\d''$ be the derivation \eqref{eq:thedifferential}. To show 
$\d'=\d''$, it suffices to prove that the two derivations agree
on $\N_{\on{hor}}\otimes 1$, and that 
$[\d',1\otimes y^a]=[\d'',1\otimes y^a]$. 

Let $x\in \N_{\on{hor}}$. Since \eqref{eq:isso} takes $x\otimes 1$ to 
$x$, $\d'(x\otimes 1)$ is the inverse image 
of $\d^\N x=\d^\N_{\on{hor}}x+\sum_a y^a. L^\N(e_a)x$. Thus
\begin{equation}\label{eq:ok} 
\d'(x\otimes 1)=\d^\N_{\on{hor}} x\otimes 1
+(-1)^{|x|}\sum_a L^\N(e_a)x\otimes y^a=\d''(x\otimes 1).
\end{equation}
We next observe that 
$ [\d^\N,y^a]=\d^W y^a=v^a+\d^\wedge y^a$
since $\N$ is a differential $W\g$-module. Hence, 
\begin{equation}\label{eq:ynot}
 [\d',1\otimes y^a]=v^a\otimes 1+1\otimes \d^\wedge y^a
=[\d'',1\otimes y^a].\end{equation}
\end{proof}

\subsection{The Chevalley-Koszul complex}
The {\em Chevalley-Koszul complex} of a $\g$-differential $W\g$-module
$\N$ is the differential $(\wedge\g)_{\on{inv}}$-module,
\begin{equation}\label{eq:chevalley}
\N_{\on{basic}}\otimes (\wedge\g^*)_{\on{inv}},\ \ \d\otimes  1+\sum_j
p^j\otimes \iota^\wedge(c_j). 
\end{equation}
\begin{remark}
If $\g$ is Abelian, \eqref{eq:chevalley} is the same as 
$(\N_{\on{hor}}\otimes\wedge\g^*)_{\on{inv}}$ with the differential \eqref{eq:thedifferential}.
\end{remark}

Let $f\in (S\g^*\otimes(\wedge\g)^-)_{\on{inv}}$ be a degree $0$
solution of the Maurer-Cartan equation \eqref{eq:basic}. Given a
$\g$-differential $W\g$-module $\N$, let $\iota^\N(f)\in \on{End}(\N)$ 
be the nilpotent operator, where the $S\g^*$-factor acts via the 
$W\g$-module structure, and the $\wedge\g$ factor acts by contraction. 
(In the special case $\N=W\g$, we write $\iota^W(f)$.)
Let $\alpha,\beta\in
\on{End}(\N_{\on{inv}}\otimes (\wedge\g^*)_{\on{inv}})$ be the
nilpotent endomorphisms of degree $0$, 
$$\alpha=\sum_j \iota(c_j)\otimes c^j,\ \beta=\sum_j c^j\otimes \iota(c_j).
$$
\begin{theorem}\label{th:chevalley}
The Chevalley-Koszul complex of $\N$ is homotopy equivalent, as a 
differential $(\wedge\g)_{\on{inv}}$-module, to the invariant 
subcomplex of $\N$. In more detail, 
\begin{enumerate}
\item\label{it:aa}
There is a homomorphism of differential $(\wedge\g)_{\on{inv}}$-modules,
$$\Psi\colon \N_{\on{basic}} \otimes (\wedge\g^*)_{\on{inv}} \to
\N_{\on{inv}},\ z\otimes\eta \mapsto
(-1)^{|\eta||z|}(e^{\iota^W(f)}\eta).z.$$
\item\label{it:bb}
There is a homomorphism of differential $(\wedge\g)_{\on{inv}}$-modules,
$$\Upsilon\colon \N_{\on{inv}}\to \N_{\on{basic}}
\otimes(\wedge\g^*)_{\on{inv}} ,\ z\mapsto (P_{\on{hor}}\otimes1
)\circ e^{-\alpha} (e^{-\iota^\N(f)}z\otimes 1).
$$
\item\label{it:cc} 
The composition $\Upsilon\circ \Psi$ is equal to the identity, while 
$\Psi\circ \Upsilon$ is homotopic to the identity by a 
$(\wedge\g)_{\on{inv}}$-homotopy. 
\end{enumerate}
\end{theorem}

\begin{proof}
\eqref{it:aa} 
Write $\iota\colon S\g^*\otimes \wedge\g\to \on{End}(\N\otimes \wedge\g^*)$, 
where the first factor acts on $\N$ via the $W\g$-module structure, and the
second factor acts by contraction.  The map $\Psi$ can be written as a
composition of two homomorphisms of $(\wedge\g)_{\on{inv}}$-modules,
$$ \N_{\on{basic}}\otimes (\wedge\g^*)_{\on{inv}}
\stackrel{e^{\iota(f)}}{\lra}
(\N_{\on{hor}}\otimes (\wedge\g^*))_{\on{inv}}
\to \N_{\on{inv}},$$
where the second 
map is an isomorphism given by \eqref{eq:canonicalisomorphism}. 
Notice that the last term in the formula \eqref{eq:thedifferential}
on $\N_{\on{hor}}\otimes (\wedge\g^*)$ vanishes on invariants. Hence, the 
differential on $(\N_{\on{hor}}\otimes (\wedge\g^*))_{\on{inv}}$ is
\[  -1\otimes \d^\wedge+\d_{\on{hor}}\otimes 1
+\sum_a v^a\otimes \iota^\wedge(e_a).\]
Conjugation of this result by the automorphism $e^{-\iota(f)}$ gives, 
\begin{equation}\label{eq:afterconjugation}
-1\otimes \d^\wedge+\d_{\on{hor}}\otimes1 +\sum_j p^j \otimes \iota^\wedge(c_j)
+\sum_a \iota(\iota^*(e^a)f)\,L^\N(e_a).
\end{equation}
On the subspace $\N_{\on{basic}}
\otimes (\wedge\g^*)_{\on{inv}}$, both 
$L^\N(e_a)$ and $1\otimes \d^\wedge$ vanish 
while $\d_{\on{hor}}\otimes 1$
coincides with $\d\otimes 1$. 
This shows that $\Psi$ is a cochain map, proving \eqref{it:aa}. 
As a preparation for \eqref{it:cc}, let us also remark that $\Psi$ is a 
homotopy equivalence of differential $(\wedge\g)_{\on{inv}}$-modules, 
by an argument parallel to the proof of Theorem \ref{th:cartan}. 
That is, starting with the operator $h= 1\otimes\sum_a
\iota(B^\sharp(e^a)) L(e_a)$ one constructs a homotopy on $(
\N_{\on{hor}}\otimes\wedge\g^*)_{\on{inv}}$ between the identity and
some projection onto $\N_{\on{basic}}\otimes(\wedge\g^*)_{\on{inv}} $.  

\eqref{it:bb} 
Consider the tensor product $ K(\P) \otimes \N_{\on{inv}}$
with differential, $\d^K\otimes 1+1\otimes \d^\N$ and with the 
$\wedge \P$-module structure 
given by $\iota(c_j)\otimes 1+1\otimes \iota(c_j)$. 
Clearly, the inclusion $z\mapsto 1\otimes z$ of $\N_{\on{inv}}$ 
is a homomorphism of differential 
$(\wedge\g)_{\on{inv}}$-modules. Under the automorphism 
$\exp(\beta)$ of $K(\P) \otimes \N_{\on{inv}}$, 
\[
\begin{split}
e^{\beta}\circ 
(\iota^K(c_k)\otimes 1+1\otimes \iota^\N(c_k))
\circ e^{-\beta}
&=\iota^K(c_k)\otimes 1,\\
e^{\beta}\circ 
(\d^K\otimes 1+1\otimes \d^\N)
\circ e^{-\beta}
&=\d^K\otimes 1+1\otimes \d^\N-\sum_j p^j\otimes \iota(c_j).
\end{split}
\]
Writing $K(\P)=(S\g^*)_{\on{inv}}\otimes (\wedge\g^*)_{\on{inv}}$, 
and re-arranging the factors, 
we hence obtain an isomorphism of $K(\P) \otimes \N_{\on{inv}}$ with 
$\ti{C}_\g(\N)\otimes (\wedge\g^*)_{\on{inv}}$, with 
differential 
$$ \ti{\d}_\g\otimes 1+\sum_j p^j\otimes \iota^\wedge(c_j), $$
and with $(\wedge\g)_{\on{inv}}$-module structure given by by
contractions on the first factor. By Theorem \ref{th:cartan} together with
Cartan's Theorem \ref{th:cargu}, there is a homotopy equivalence of
differential $(S\g^*)_{\on{inv}}$-modules,
$$ (S\g^*)_{\on{inv}}\otimes \N_{\on{inv}}
=\ti{C}_\g(\N)\to C_\g(\N) \to \N_{\on{basic}},\ \ 
p\otimes z\mapsto p.P_{\on{hor}}(e^{\iota^W(f)}.z).
$$
The resulting homomorphism of differential $(\wedge\g)_{\on{inv}}$-modules, 
$$ \N_{\on{inv}}\to 
(\wedge\g^*)_{\on{inv}}\otimes \big((S\g^*)_{\on{inv}}\otimes \N_{\on{inv}}\big)\to 
(\wedge\g^*)_{\on{inv}}\otimes \N_{\on{basic}}$$
is exactly our map $\Upsilon$. 

\eqref{it:cc} 
Let $z\otimes\eta
\in  \N_{\on{basic}}\otimes (\wedge\g^*)_{\on{inv}}$. We compute, 
\[
\begin{split}
\Upsilon(\Psi(z\otimes\eta ))&=(-1)^{|z||\eta|}
\Upsilon\big((e^{\iota^W(f)}\eta).z\big)\\
&=
(-1)^{|z||\eta|}(P_{\on{hor}}\otimes 1 )\circ 
e^{-\alpha}e^{-\iota^{\N}(f) \otimes 1} ( (e^{\iota^W(f)}\eta).z\otimes 1)\\
&=(-1)^{|z||\eta|}(P_{\on{hor}}\otimes 1 )\circ 
e^{-\alpha}
( \eta.z \otimes 1)\\
&=z\otimes \eta.
\end{split}\]
Thus $\Upsilon\circ \Psi=\on{Id}$, while the opposite composition
$\Pi=\Psi\circ\Upsilon$ is a projection.  As remarked in \eqref{it:aa} above,
there exists a $(\wedge\g)_{\on{inv}}$-homotopy operator $H_1$ between
$I$ and {\em some} projection operator $\Pi_1$ onto
$\N_{\on{basic}}\otimes(\wedge\g^*)_{\on{inv}}$.  Lemma
\ref{lem:easy3} below shows how to obtain from this a
$(\wedge\g)_{\on{inv}}$-homotopy operator $H$ between $I$ and $\Pi$.
\end{proof}

\begin{lemma}\label{lem:easy3}
Let $(C,\d)$ be an $\A$-differential space, where $\A$ is 
some graded algebra (with trivial differential). 
Suppose $H_1\colon C\to C$ is an $\A$-equivariant homotopy 
operator between the identity and some $\A$-equivariant 
projection $\Pi_1\colon C\to C$ onto a differential 
subspace $C'\subset C$.  Assume that some cochain map 
$\Pi\colon C\to C$ is another $\A$-equivariant projection 
onto $C'$. Then 
$$ H=(I-\Pi_1-\Pi)H_1(I-\Pi_1-\Pi)$$ 
is an $\A$-equivariant homotopy between $I$ and $\Pi$. 
\end{lemma}
\begin{proof}
This follows by straightforward calculation, using $\Pi_1\Pi=\Pi$ and 
$\Pi\,\Pi_1=\Pi_1$:
\[\begin{split}
 [\d,H]&=(I-\Pi_1-\Pi)[\d,H_1](I-\Pi_1-\Pi)\\
&=(I-\Pi_1-\Pi)(I-\Pi_1)(I-\Pi_1-\Pi)=I-\Pi.\end{split}\]
\end{proof}

\begin{remark}
Suppose $f_0,f_1$ are two solutions of \eqref{eq:basic}.
Then, the corresponding cochain maps from 
$(\wedge\g^*)_{\on{inv}}\otimes\N_{\on{basic}}$ into $\N_{\on{inv}}$ are 
homotopic,  by a homotopy operator which is compatible with the 
$(\wedge\g^*)_{\on{inv}}$-module structure. The proof is parallel to that of 
Proposition \ref{prop:homotopy}. 
\end{remark}

\begin{remark}
Recall that in \ref{rem:chevhom}, we described a homomorphism of
differential graded algebras $\Phi'\colon K(P)\to (W\g)_{\on{inv}}$, 
depending on the choice of distinguished cochains of transgression. Thinking of 
$(\wedge\g^*)_{\on{inv}}=\wedge P^*$ as a subalgebra of $K(P)$, this gives a
cochain map, 
\[ \Psi'\colon \N_{\on{basic}} \otimes (\wedge\g^*)_{\on{inv}}\stackrel{1\otimes \Phi'}{\lra}
 \N_{\on{basic}} \otimes (W\g)_{\on{inv}}
\to \N_{\on{inv}},\]
where the last map is $x\otimes w\mapsto (-1)^{|x||w|}w.x$. This map
is described in \cite[page 364]{gr:co3} under the name \emph{Chevalley
homomorphism}. Assuming that the complex $\N$ is bounded below, 
an easy spectral sequence argument shows that $\Psi'$
induces an isomorphism in cohomology \cite[page 365]{gr:co3}. In
\cite{mas:ko}, it is claimed that $\Psi'$ is a homomorphism of $\wedge
P$-modules, but this is false (see Remark \ref{rem:chevhom}).

The fact that the complex $\N_{\on{basic}}\otimes(\wedge\g^*)_{\on{inv}}$ 
computes the cohomology
of $\N_{\on{inv}}$ (hence of $\N$, since $\g$ is reductive) goes back
to Chevalley and Koszul.  See the article of Koszul \cite{ko:tr1} in
the `Colloque de topologie' (reproduced in \cite{ko:se}).  As a special case of
this result, one obtains the de Rham cohomology of the total space of
a principal $G$-bundle $P\to B$ (with $G$ a compact Lie group) as
the cohomology of a complex $\Om(B)\otimes (\wedge\g^*)_{\on{inv}}$,
where the complex $\Om(B)$ of differential forms on the base is viewed
as an $(S\g^*)_{\on{inv}}$-module by the Chern-Weil homomorphism.
\end{remark}

As explained in \cite{gor:eq,mas:ko}, Theorems \ref{th:chevalley} and
\ref{th:cartan} are related by Koszul duality. We outline this argument
in Appendix B. In Appendix A we indicate a common generalization of 
the two theorems.

\section{Lie algebra homomorphisms}\label{sec:liehom}
In this Section, we will address the functoriality properties of the
small Cartan model and of the Chevalley-Koszul model under
homomorphisms of reductive Lie algebras. In particular, taking $\g\to
\g\oplus\g$ the diagonal embedding, this will also lead to a
description of product structures for the two models.
\subsection{Lie algebra homomorphisms}
Suppose $\h$ is another reductive Lie algebra, with $\ca{Q}\subset
(\wedge\h)_{\h\mbox{-}\!\on{inv}}$ as its space of primitive elements, and that
$\phi\colon\g\to \h$ is a Lie algebra homomorphism.  Then $\phi$ extends
to an algebra homomorphism $\phi\colon\wedge\g\to \wedge\h$, compatible
with the boundary operator and with the Schouten bracket. Furthermore, 
the dual map restricts to a linear map 
$$ \phi^*\colon(\wedge\h^*)_{\h\mbox{-}\!\on{inv}}\to (\wedge\g^*)_{\g\mbox{-}\!\on{inv}},$$
which is a homomorphism of Hopf algebras \cite[Section 5.17]{gr:co3}.  
Hence it restricts to a linear map, 
$$ \phi^*\colon  \ca{Q}^*\to  \P^*.$$
Let $d_l\in \ca{Q}$ and $q^l\in \ti{\ca{Q}}^*$ by dual 
(homogeneous) bases.

\begin{proposition}\label{prop:fun}
There exists a degree $0$ solution $u\in (S\g^*)_{\g\mbox{-}\!\on{inv}}\otimes 
(\wedge\h)^-_{\g\mbox{-}\!\on{inv}}$ of the equation, 
\begin{equation}\label{eq:u}
 \partial u +\hh [u,u]_{\wedge\h}=
\sum_l \phi^*(q^l) \otimes d_l-\sum_j p^j \otimes \phi(c_j).
\end{equation}
\end{proposition}
\begin{proof}
We will apply Theorem \ref{th:general} to the setting, 
\[\begin{split}
\k&=\bigoplus_{i\le 0} \k^i,\ \ \k^i=(S\g^*)_{\g\mbox{-}\!\on{inv}}\otimes (\wedge^{1-i}\h)_{\g\mbox{-}\!\on{inv}},\\
\mf{l}&=\bigoplus_{i\le 0}\mf{l}^i,\ \   
\mf{l}^i= (S\g^*)_{\g\mbox{-}\!\on{inv}}\otimes 
(\wedge^{1-i}\h)_{\h\mbox{-}\!\on{inv}},\\
X&=\sum_l \phi^*(q^l) \otimes d_l-\sum_j p^j \otimes \phi(c_j).
\end{split}\]
Here the bracket $[\cdot,\cdot]_\k$ and differential $\partial$ on
$\k$ are induced from the Schouten bracket and differential on
$(\wedge\h)_{\g\mbox{-}\on{inv}}\subset \wedge\h$.  Observe that
$\mf{l}$ is contained in the center of $\k$, since 
$(\wedge\h)_{\h\mbox{-}\on{inv}}$ is the center of $\wedge\h$.  
The differential vanishes on $\mf{l}$, and by Hodge theory for $\wedge\h$ the inclusion
$\mf{l}\hra \k$ induces an isomorphism in cohomology. Finally $X$ is a
cocycle (since $\partial q^l=0$ and $\partial c^j=0$), and is contained in
the center since elements in $\phi(\wedge\g)$ Schouten commute with
$\g$-invariant elements in $\wedge\h$.  Hence all assumptions of
Theorem \ref{th:general} are satisfied, and we obtain a solution $u$
of \eqref{eq:u} {\em modulo $\mf{l}$}, (Using Hodge theory on
$\wedge\h$, we obtain in fact a {\em canonical} solution $u$ of
(total) degree $0$.) Theorem \ref{th:general} says furthermore that
the `error term'
$$Y:=\partial u +\hh [u,u]_{\wedge\g}+
\sum_l \phi^*(q^l) \otimes d_l-\sum_j p^j \otimes \phi(c_j)
\in\mf{l}$$
does not depend on the choice of $u$.  It remains to show that in 
fact $Y=0$.

Let $f\in (S\g^*\otimes(\wedge\g)^-)_{\g\mbox{-}\!\on{inv}}$
be an even element solving \eqref{eq:basic}. Using again that 
elements in $\phi(\wedge\g)$ and $(\wedge\h)_{\g\mbox{-}\!\on{inv}}$
commute under the Schouten bracket, we have
\[ [(1\otimes\phi)(f),u]_{\wedge\h}=0. \]
It follows that
\[\ti{f}_1:=(1\otimes \phi)(f)+u\in (S\g^*\otimes(\wedge\h)^-)_{\g\mbox{-}\!\on{inv}}\]
satisfies the equation 
\begin{equation}\label{eq:theqn}
 \partial \ti{f}_1+\hh [\ti{f}_1,\ti{f}_1]_{\wedge\h}=
\sum_l \phi^*(q^l)\otimes  d_l-\sum_a v^a \otimes \phi(e_a)+Y.
\end{equation}
On the other hand, let 
$g\in (S\h^*\otimes(\wedge\h)^-)_{\h\mbox{-}\!\on{inv}}$ be 
an even solution of the analogue of \eqref{eq:basic} for the 
Lie algebra $\h$. Then $\ti{f}_0=(\phi^*\otimes 1)(g)$ 
solves a similar equation, but with $Y$ replaced by $0$:
\begin{equation}\label{eq:theqn1}
 \partial \ti{f}_0+\hh [\ti{f}_0,\ti{f}_0]_{\wedge\h}=
\sum_l \phi^*(q^l)\otimes  d_l-\sum_a v^a \otimes \phi(e_a).
\end{equation}
(Here we used that the elementary fact that the image of the 
canonical element $\sum_a v^a \otimes e_a\in S\g^*\otimes \wedge\g$ 
under $1\otimes \phi$ coincides with the image of the corresponding 
element of $S\h^*\otimes\wedge\h$ under the map $\phi^*\otimes 1$.) 
By the uniqueness part of Theorem \ref{th:general}, applied to the situation
\[\begin{split}
\k&=\bigoplus_{i\le 0}\k^i,\ \  
\k^i=(S\g^*\otimes (\wedge\h))_{\h\mbox{-}\!\on{inv}},\\
\mf{l}&=\bigoplus_{i\le 0}\mf{l}^i,\ \
\mf{l}^i= (S\g^*)_{\g\mbox{-}\!\on{inv}}\otimes 
(\wedge^{1-i}\h)_{\h\mbox{-}\!\on{inv}},\\
X&=\sum_l \phi^*(q^l)\otimes  d_l-\sum_a v^a \otimes \phi(e_a),
\end{split}\]
this shows $Y=0$. 
\end{proof}

\begin{remark}
If the Lie algebra $\g$ is {\em Abelian}, one may simply take
$u=(\phi^*\otimes 1)(g)$ where $g\in
(S\h^*\otimes(\wedge\h)^-)_{\h\mbox{-}\on{inv}}$ solves
\eqref{eq:basic2}.  For instance, let $\h$ be the 3-dimensional Lie
algebra from Example \ref{ex:su2} (denoted $\g$ in that example),
and let $\g\subset \h$ be the inclusion of the 1-dimensional Lie subalgebra
spanned by $e_3$. Then $u=v^3 \otimes (e_1\wedge e_2)$ is a solution
of \eqref{eq:u}.
\end{remark}

As a special case of Proposition \ref{prop:fun}, 
consider the diagonal inclusion $\on{diag}\colon \g\to
\g\oplus\g$, and its extension to the exterior algebra,
$\on{diag}\colon \wedge\g\to \wedge\g\otimes\wedge\g$. The dual map
$\on{diag}^*$ is just the product map for $\wedge\g^*$. Hence 
Proposition \ref{prop:fun}
defines a solution $u\in (S\g^*)_{\g\mbox{-}\!\on{inv}}\otimes
(\wedge\g\otimes\wedge\g)_{\g\mbox{-}\!\on{inv}}$ of the equation,
\begin{equation}\label{eq:product}
\partial u+\hh[u,u]_{\wedge(\g\oplus\g)}=
\sum_j p^j \otimes (\on{diag}(c_j)-\Delta(c_j))
\end{equation}
where $\Delta$ 
is the coproduct on $(\wedge\g)_{\on{inv}}$, i.e. 
$\Delta(c_j)=c_j\otimes 1+1\otimes c_j$. 

\begin{example}
Consider $\g$ as in Example \ref{ex:su2}, with 
$p=\sum_a v^a v^a$.  Write $e_a^{(1)}$ (resp. $e_a^{(2)}$) for the 
basis vectors in the first (resp. second) copy of $\g$ in $\g\oplus\g$. 
Then 
$$ u=p\otimes (e_1^{(1)}\wedge e_2^{(1)}\wedge e_1^{(2)}\wedge e_2^{(2)}
+\ldots)$$
(where the dots indicate a sum over cyclic permutations over 
the lower indices $1,2,3$) solves Equation \eqref{eq:product}. Note that 
$[u,u]_{\wedge(\g\oplus\g)}=0$ in this case. 
\end{example}

\subsection{Small Cartan complex}
Let $\phi\colon\g\to \h$ is a homomorphism of reductive Lie 
algebras, and $\M$ be an $\h$-differential space. The natural map 
$$ C_\h(\M)=(S\h^*\otimes \M)_{\h\mbox{-}\!\on{inv}}\to 
C_\g(\M)=(S\g^*\otimes\M)_{\g\mbox{-}\!\on{inv}}$$
is a cochain map, inducing a map in cohomology, 
$H_\h(\M)\to H_\g(\M)$. We will now realize this map in terms of 
the small Cartan models. 
\begin{theorem}\label{th:fun1}
Let $u\in (S\g^*)_{\g\mbox{-}\!\on{inv}}\otimes 
(\wedge\h)_{\g\mbox{-}\!\on{inv}}$ be 
a solution of \eqref{eq:u}. Then the operator 
$$ \Psi=e^{\iota(u)}\circ (\phi^*\otimes 1)\colon (S\h^*)_{\h\mbox{-}\!\on{inv}}
\otimes \M_{\h\mbox{-}\!\on{inv}}\to (S\g^*)_{\g\mbox{-}\!\on{inv}}
\otimes \M_{\g\mbox{-}\!\on{inv}}
$$
is a cochain map. The resulting diagram of differential 
$(S\g^*)_{\g\mbox{-}\!\on{inv}}$-modules, 
\begin{equation}
\label{eq:diagram}
\xymatrix{ {(S\h^*)_{\h\mbox{-}\!\on{inv} }\otimes
\M_{\h\mbox{-}\!\on{inv} } } \ar[r]\ar[d] & { (S\h^*\otimes
\M)_{\h\mbox{-}\!\on{inv} } } \ar[d] \\ {
(S\g^*)_{\g\mbox{-}\!\on{inv} }\otimes \M_{\g\mbox{-}\!\on{inv} }
}\ar[r]& { (S\g^*\otimes \M)_{\g\mbox{-}\!\on{inv} } } }
\end{equation}
commutes up to $(S\g^*)_{\g\mbox{-}\!\on{inv}}$-equivariant homotopy. In particular,
the map $\ti{H}_\h(\M)\to \ti{H}_\g(\M)$ defined by $\Psi$ agrees with the 
natural map $H_\h(\M)\to H_\g(\M)$, under the identification 
$\ti{H}_\g(\M)=H_\g(\M)$, $ \ti{H}_\h(\M)=H_\h(\M)$ from Theorem 
\ref{th:cartan}. 
\end{theorem}
\begin{proof}
Since the image of $(S\h^*)_{\h\mbox{-}\!\on{inv}}\otimes \M_{\h\mbox{-}\!\on{inv}}$ under the map $(\phi^*\otimes 1)$ is contained in
$(S\g^*)_{\g\mbox{-}\!\on{inv}}\otimes \M_{\h\mbox{-}\!\on{inv}}$,
\eqref{eq:comm2} and 
\eqref{eq:u} show that
$$(1\otimes\d)\circ \Psi= \Psi \circ \big(1\otimes\d -
\sum_l \phi^*(q^l)\otimes  d_l-\sum_a v^a \otimes \phi(e_a)
\big).$$
Hence $\Psi$ is a cochain map. Now let $f\in
(S\g^*\otimes(\wedge\g)^-)_{\on{inv}}$ be an even solution of \eqref{eq:basic},
and $g\in (S\h^*\otimes (\wedge\h)^-)_{\h\mbox{-}\!\on{inv}}$ an even solution
of the corresponding equation for $\h$. 

As we observed in the proof of Proposition \ref{prop:fun}, the two elements
$$ \ti{f}_0=(\phi^*\otimes 1)g,\ 
\ti{f}_1=(1\otimes \phi)(f)+u$$
of $(S\g^*\otimes(\wedge\h)^-)_{\g\mbox{-}\!\on{inv}}$
both satisfy the equation
$$ 
 \partial \ti{f}+\hh [\ti{f},\ti{f}]_{\wedge\h}=
\sum_l \phi^*(q^l)\otimes  d_l-\sum_a v^a \otimes \phi(e_a).
$$
By Theorem \ref{th:general}, it follows that 
$\ti{f_0}$ and $\ti{f}_1$ are gauge
equivalent, up to addition of an even element in 
$(S\g^*)_{\g\mbox{-}\!\on{inv}}\otimes 
(\wedge^{1-i}\h)_{\h\mbox{-}\!\on{inv}}$.
The argument at the beginning of Section \ref{sec:choicef} shows
that the action by elements in $(S\g^*)_{\g\mbox{-}\!\on{inv}}\otimes 
(\wedge^{1-i}\h)_{\h\mbox{-}\!\on{inv}}$ is trivial, up to homotopy. 
Hence we can restrict consideration to gauge equivalent
$\ti{f}_1=\exp(s_1).\ti{f}_0$, where 
$s\in (S\g^*\otimes(\wedge\h)^-)_{\g\mbox{-}\!\on{inv}}$ is odd. 
These then belong to a family 
$\ti{f}(t)=\exp(s(t)).\ti{f}_0$ of solutions, where 
$s(t)=t s_1$. We have, 
$$ \f{d \ti{f}}{d t}+\partial s+[\ti{f},s]_{\wedge\h}=0.$$
Each $\ti{f}(t)$ defines a cochain map
$$ \ti{\Psi}(t)=e^{\iota(\ti{f}(t))}\circ (\phi^*\otimes 1)\colon 
{(S\h^*)_{\h\mbox{-}\!\on{inv}}\otimes \M_{\h\mbox{-}\!\on{inv}}}\to 
(S\g^*\otimes \M)_{\g\mbox{-}\!\on{inv}}.$$
Consider the family of operators $ H(t)=\ti{\Psi}(t)\circ \iota(s(t))$. 
Arguing as in the proof of Proposition \ref{prop:homotopy}, one shows that 
$$ \f{d \ti{\Psi}}{d t}=H\circ \ti{\d}_\h+\d_\g\circ H.$$ 
Hence the maps $\ti{\Psi}(t)$ are all homotopic, where the homotopy 
respects the $(S\g^*)_{\g\mbox{-}\!\on{inv}}$-module structures.  
\end{proof}

Suppose now that $\M$ is a $\g$-differential {\em algebra}, that is
$\M$ is a $\g$-differential space which is also a graded algebra, 
in such a way that $\d,\iota(\xi),L(\xi)$ are all derivations for the product 
$\mu_\M\colon \M\otimes\M\to\M$. In this case, the Cartan complex $C_\g(\M)$ 
inherits a product structure for which $\d_\g$ is a derivation. 
Hence, the product descends to the cohomology $H_\g(\M)$.

By contrast, the differential $\ti{\d}_\g$ for the small Cartan model is
not a derivation for the obvious product structure on
$(S\g^*)_{\on{inv}}\otimes \M_{\on{inv}}$.  Instead, define a new
(non-associative) multiplication $\odot$ on 
$(S\g^*)_{\on{inv}}\otimes \M_{\on{inv}}$ by 
\[\begin{split}
(p\otimes y)\odot (p'\otimes y')=
(1\otimes \mu_\M) e^{\iota(u)}(pp'\otimes y\otimes y'),
\end{split}
\]
where $u\in (S\g^*)_{\g\mbox{-}\!\on{inv}}\otimes
(\wedge\g\otimes\wedge\g)_{\g\mbox{-}\!\on{inv}}$ is a solution of 
\eqref{eq:product}.

\begin{theorem}\label{th:products}
The differential $\ti{\d}_\g$ on the small Cartan model is a derivation with 
respect to the product $\odot$, 
$$ \ti{\d}_\g(x\odot x')=\ti{\d}_\g(x)\odot x'+
(-1)^{|x|}x\odot \ti{\d}_\g(x').$$
The induced product on the equivariant cohomology of $\M$ coincides with that 
from the usual Cartan model.  
\end{theorem}
\begin{proof}
This follows from Theorem \ref{th:fun1}, specialized 
to the diagonal inclusion $\g\to \g\oplus\g$. 
\end{proof}

\subsection{Chevalley-Koszul complex}
Suppose $\phi\colon \g\to\h$ is a homomorphism of reductive Lie algebras. 
The dual map 
$\phi^*\colon \h^*\to\g^*$ extends to a homomorphism of 
$\g$-differential algebras
$$\phi^*\colon W\h\to W\g.$$
Suppose $\M$ is an $\g$-differential $W\g$-module, $\N$ is an
$\h$-differential $W\h$-module, and $F\colon \N\to \M$ is a homomorphism of 
$\g$-differential spaces such that the following diagram commutes:
$$ \xymatrix{W\h\otimes \N\ar[r]\ar[d]_{\phi^*\otimes F}
&\N\ar[d]^F\\ W\g\otimes \M\ar[r] &\M}$$
The geometric setting to have in mind is that of a reduction of the structure 
group of a principal with connection. 
The following result may be proved along the same lines as 
Theorem \ref{th:fun1}:
\begin{theorem}\label{th:fun2}
Let $u\in (S\g^*)_{\g\mbox{-}\!\on{inv}}\otimes 
(\wedge\h)_{\g\mbox{-}\!\on{inv}}$ be a solution of 
\eqref{eq:u}. Then the composition of maps
$$
(F\otimes 1)\circ e^{\iota(u)}\circ (1\otimes \phi^*)\colon 
 \N_{\h\mbox{-}\!\on{basic}} \otimes (\wedge\h^*)_{\h\mbox{-}\!\on{inv}}
\to \M_{\g\mbox{-}\!\on{basic}}\otimes (\wedge\g^*)_{\g\mbox{-}\!\on{inv}}$$
is a homomorphism of differential
$(\wedge\g)_{\h\mbox{-}\!{\on{inv}}}$-modules.  Under the isomorphism
from Theorem \ref{th:chevalley}, this map induces the same map in
cohomology as the cochain map $F\colon\N_{\h\mbox{-}\!\on{inv}}\to
\M_{\g\mbox{-}\!\on{inv}}$.
\end{theorem}

\begin{appendix}
\section{The Halperin complex}\label{sec:halperin}
Suppose $P\to B$ is a principal $G$-bundle, $F$ any $G$-manifold, and
$P\times_G F$ the associated bundle with fiber $F$. An unpublished
result of Halperin (quoted in \cite[page 569]{gr:co3}) describes the
de Rham cohomology of $P\times_G F$ as the cohomology of a certain differential 
on $\Om(B)\otimes \Om(F)_{\on{inv}}$. 
The following is an algebraic generalization of Halperin's
result. 

\begin{theorem}\label{th:halperin}
Suppose that $\N$ is a $\g$-differential $W\g$-module, 
and that $\M$ is any $\g$-differential space. Let 
$\N_{\on{basic}}\otimes \M_{\on{inv}}$ 
be equipped with the differential, 
\begin{equation}\label{eq:halperin}
 \d^\N\otimes 1+1\otimes \d^\M-\sum p^j\otimes \iota^\M(c_j).
\end{equation}
Define a nilpotent operator of degree 0, 
$\alpha=\sum_a y^a\otimes\iota^\M(e_a)\in \on{End}(\N\otimes\M)$, 
and let $f\in (S\g^*\otimes(\wedge\g)^-)_{\on{inv}}$ be
a degree $0$ element solving \eqref{eq:basic}.
Then the map 
\begin{equation}\label{eq:halp1}
e^\alpha \circ e^{\iota^\M(f)}\colon  
\N_{\on{basic}}\otimes \M_{\on{inv}}\to (\N\otimes\M)_{\on{basic}},
\end{equation}
is a cochain map. It is a homotopy equivalence of differential
$(S\g^*)_{\on{inv}}$-modules, provided at least one of the
$\g$-differential spaces $\N,\ \M$ is a direct sum of the kernel and
image for the action of the Casimir operator $\on{Cas}_\g$.
\end{theorem}

\begin{proof}
The operator $e^{-\alpha}$ on $\N\otimes\M$
takes the contractions on $\N\times\M$ to contractions on the first factor, 
$$e^\alpha\circ (\iota^\N(\xi)+\iota^\M(\xi))\circ e^{-\alpha}= 
\iota^\N(\xi).$$
In particular, it restricts to an isomorphism $ (\N\otimes\M)_{\on{basic}}\to
(\N_{\on{hor}}\otimes \M)_{\on{inv}}.$
For the induced differential on
$(\N_{\on{hor}}\otimes \M)_{\on{inv}}$ one finds, after short
calculation (cf. \cite[page 18]{ka:br})
$$ \d':=
e^\alpha\circ (\d^\N+ \d^\M)\circ 
e^{-\alpha}=
\d^\N_{\on{hor}}+\d^\M
-\sum_a v^a\,\iota^{\M}(e_a).$$
(If $\N=W\g$, this is Kalkman's proof \cite{ka:br} of the
equivalence between $(W\g\otimes\M)_{\on{basic}}$ and the Cartan
model $(S\g^*\otimes\M)_{\on{inv}}$.)
The operator
$\exp(\iota^\M(f))$ commutes with
$\iota^\N(\xi)$, and therefore preserves $(\N_{\on{hor}}\otimes
\M)_{\on{inv}}$. It also commutes with $\d^\N_{\on{hor}}$ and 
$\sum_a v^a\,\iota^{\M}(e_a)$, hence the only new contributions arise from
the commutator with $\d^\M$. We find,
\begin{equation}
\d'':=e^{-\iota^\M(f)}\circ \d'\circ
e^{\iota^\M(f)}
=\d^\N_{\on{hor}}+\d^\M
-\sum_j p^j \iota^{\M}(c_j)+\sum_a \iota^\M(\iota^*(e^a) f)\circ 
L^\M(e_a).
\end{equation}
On the subalgebra $\N_{\on{basic}}\otimes \M_{\on{inv}}\subset
(\N_{\on{hor}}\otimes \M)_{\on{inv}}$, the differential simplifies to
\eqref{eq:halperin}, as desired. Working backwards, this defines the
cochain map \eqref{eq:halp1}. To show that \eqref{eq:halp1} is a
homotopy equivalence, it suffices to show that the inclusion of
$\N_{\on{basic}}\otimes \M_{\on{inv}}$ into 
$(\N_{\on{hor}}\otimes \M)_{\on{inv}}$ is a homotopy equivalence.
This is done by a straightforward extension of the argument given in the 
proof of Theorem \ref{th:cartan}. The only fact needed is that 
$(\N_{\on{hor}}\otimes \M)_{\on{inv}}$ is a direct sum of the kernel and 
image of $\ca{L}_0=\on{Cas}_\g^\N\otimes 1$, which follows by assumption. 
(Note that  $\ca{L}_0=1\otimes \on{Cas}_\g^\M$ on invariants.)
\end{proof}

Note that Theorem \ref{th:halperin} specializes to Theorem 
\ref{th:cartan} for $\N=W\g$. 
Similarly, it contains 
part of Theorem \ref{th:chevalley} as the special case 
$\M=\wedge\g^*$.

\section{Koszul duality}\label{app:koszul}
In this appendix, we discuss the Koszul duality between differential
graded modules over symmetric and exterior algebras. Much of this
discussion is already implicit in Koszul's work \cite{ko:tr1}, and has
appeared in the literature in various degrees of generality
\cite{gor:eq,gr:co3,hue:ho,hue:re}. For simplicity, all differential
complexes in this appendix are assumed to be bounded below,
i.e. equal to 0 in sufficiently negative degrees. 

Let $\P$ be any finite-dimensional graded vector space, concentrated 
in odd negative degrees, and let $\P^*$ be its dual space with grading 
$(\P^*)^i=(\P^{-i})^*$. Write $\ti{\P}=\P[-1]$ and $\P^*=\P^*[1]$. 

Fix dual bases $c_j,c^j$ of
$\P,\P^*$, and let $p_j,p^j$ denote the corresponding dual bases of
$\ti{\P},\ti{\P}^*$. The Koszul complex of $\P$ is the 
differential graded algebra 
\[ K(\P)=S\ti{\P}^*\otimes \wedge \P^*,\ \ 
\d_{K(\P)}=\sum_j p^j\otimes\iota(c_j).\]
As is well-known, $K(\P)$ is acyclic.  Let $K(\P)^-$ denote the Koszul complex with differential
$\d_{K(\P)^-}=-\d_{K(\P)}$.

A {\em differential $S\ti{\P}^*$-module} is a differential space
$(X,\d_X)$, which is also an $S\ti{\P}^*$-module, in such a way that
the action of any element in $S\ti{\P}^*$ commutes with the differential.  Abusing
notation, we will denote the action of $p^j\in S\ti{\P}^*$ on $X$ simply
by $p^j$.  If $X_1,X_2$ are two differential $S\ti{\P}^*$-modules, their
tensor product $X_1\otimes X_2$ is again a differential $S\ti{\P}^*$-module. 
 
In a similar way, we define {\em differential $\wedge\P$-modules} 
$(Y,\d_Y)$. We will denote the action of $c_j \in \P$ by 
$\iota(c_j)$. Given two differential 
$\wedge\P$-modules, their tensor 
product is again a differential $\wedge\P$-module.

Note that $K(\P)$ is both a differential $S\ti{\P}^*$-module 
and a differential
$\wedge\P$-module. The augmentation map $K(\P)\to \F$ is a
homomorphism of differential $S\ti{\P}^*$-modules, while the
coaugmentation map $\F\to K(\P)$ is a homomorphism of
differential $\wedge \P$-modules. 

There is a covariant functor $h$ from the category of 
differential $S\ti{\P}^*$-modules 
to the category of differential $\wedge \P$-modules, 
taking $X$ to 
$$ hX=X\otimes \wedge \P^*,\ \ \d_{hX}=\d_X \otimes 1 +
\sum_j p^j\otimes  \iota(c_j),$$
and taking a morphism $\phi\colon X\to X'$ to $h\phi=\phi\otimes 1 \colon hX\to hX'$. 
The filtration of $hX$ coming from the grading on $\wedge\P^*$ 
defines a spectral sequence for $hX$, with $E_2$-term $H(X)\otimes \wedge\P^*$.
In particular, the functor $h$ preserves quasi-isomorphisms. 

Similarly, there is a covariant  functor $t$ from the category of differential $\wedge \P$-modules
to the category of differential $S\ti{\P}^*$-modules, taking $Y$ to  
$$ tY=S\ti{\P}^*\otimes Y,\ \ \d_{tY}=1\otimes \d_Y-\sum_j p^j\otimes
\iota(c_j).$$
and taking a morphism $\psi\colon Y\to Y'$ to $t\psi=1\otimes
\psi\colon tY\to tY'$.
Again, $t$ preserves quasi-isomorphisms. 
\begin{theorem}[Koszul duality]
\begin{enumerate}
\item
For any differential $\wedge\P$-module $Y$ 
there is a canonical isomorphism of differential $\wedge\P$-modules, 
$$ htY\cong K(\P)\otimes Y,$$
and hence a canonical quasi-isomorphism 
$ Y\to htY$. 
\item
For any differential $S\ti{\P}^*$-module $X$, there 
is a canonical isomorphism of differential $S\ti{\P}^*$-modules,
$$ thX \cong K(\P)^-\otimes X$$
and hence a canonical quasi-isomorphism $thX \to X$.
\end{enumerate} 
\end{theorem}
\begin{proof}
By definition, the differential $\wedge\P$-module 
$htY$ is equal to $K(\P)\otimes Y$ as a vector space, 
but with $\wedge\P$-module structure and differential given by 
\[ 
\begin{split}
\iota_{htY}(c_j)&=\iota_{K(\P)}(c_j)\otimes 1,\\
\d_{htY}&=\d_{K(\P)}\otimes 1+1\otimes\d_Y
-\sum_j p^j\otimes \iota_Y(c_j).
\end{split}
\]
The endomorphism $\alpha=\sum_k c^k\otimes \iota_Y(c_k)$ 
of $K(\P)\otimes Y$ is nilpotent and has degree $0$. Hence, its exponential 
is a well-defined automorphism of degree $0$. A straightforward calculation, 
using $\Ad(e^{-\alpha})=e^{-\ad(\alpha)}$, shows
\[
\begin{split}
\Ad(e^{-\alpha})(\iota_{htY}(c_j))&=\iota_{K(\P)}(c_j)\otimes 1
+1\otimes \iota_Y(c_j),\\
\Ad(e^{-\alpha})(\d_{htY})&=\d_{K(\P)}\otimes 1+1\otimes\d_Y.
\end{split}
\]
It follows that $e^\alpha\colon K(\P)\otimes Y\to htY$ is an
isomorphism of differential $\wedge\P$-modules. Similarly, $thX=
K(\P)\otimes X$ as a vector space, but with $S\ti{\P}^*$-module
structure and differential given by
\[
\begin{split}
\iota_{thX}(p^j)&=p^j\otimes 1,\\
\d_{thX}&=-\d_{K(\P)}\otimes 1+1\otimes\d_X
-\sum_j \iota_{K(Y)}(c_j) \otimes p^j.
\end{split}
\]
The endomorphism $\beta=\sum_j \iota_S(p_j)\otimes p^j$ of $K(\P)\otimes X$
is nilpotent and has degree $0$. We find, 
\[
\begin{split}
\Ad(e^\beta)\iota_{thX}(p^j)&=p^j\otimes 1+1\otimes p^j,\\
\Ad(e^\beta)\d_{thX}&=-\d_{K(\P)}\otimes 1+1\otimes\d_X.
\end{split}
\]
Hence $e^{-\beta}$ gives an isomorphism of differential 
$S\ti{\P}^*$-modules, $K(\P)^-\otimes X\to thX$.
\end{proof}

We are now in position to explain the Koszul duality between Theorems
\ref{th:cartan} and \ref{th:chevalley}. Suppose $\M$ is any 
$\g$-differential space. As before, we denote by $\P$ the primitive 
subspace of $(\wedge\g)_{\on{inv}}$. 
Then the subspace $\M_{\on{inv}}$ of invariants 
is a differential $\wedge \P$-module, while the Cartan model
$(S\g^*\otimes \M)_{\on{inv}}$ is a differential 
$S\ti{\P}^*$-module. 
Applying the functor $t$ to $\M_{\on{inv}}$, we obtain the small Cartan 
model, 
\begin{equation}\label{eq:1}
 t\M_{\on{inv}}=S\ti{\P}^*\otimes \M_{\on{inv}}.
\end{equation}
On the other hand, recall \cite{ca:la} that the map
$(W\g\otimes\M)_{\on{basic}}\to (S\g^*\otimes \M)_{\on{inv}}=C_\g(\M)$
induced by the projection $W\g\to S\g^*$ is an isomorphism 
of differential $S\ti{\P}$-modules. Hence, applying the functor 
$h$ to the (big) Cartan model we obtain 
\begin{equation}\label{eq:2}
 h(S\g^*\otimes \M)_{\on{inv}}
=  (W\g\otimes \M)_{\on{basic}}\otimes \wedge\P^*.
\end{equation}
Let $\sim$ denote the relation of quasi-isomorphism in the category of 
differential $S\ti{\P}^*$-modules, respectively of 
differential $\wedge\P$-modules. Since 
$W\g$ is acyclic, $(W\g\otimes \M)_{\on{inv}}\sim \M_{\on{inv}}.$
The quasi-isomorphism from Theorem \ref{th:chevalley}, 
$$ h(S\g^*\otimes\M)_{\on{inv}}=
h(W\g\otimes \M)_{\on{basic}}\sim 
(W\g\otimes \M)_{\on{inv}}\sim \M_{\on{inv}}$$
implies, by Koszul duality, a quasi-isomorphism of differential 
$S\ti{\P}^*$-modules, 
$$ (S\g^*\otimes\M)_{\on{inv}}\sim 
th(S\g^*\otimes\M)_{\on{inv}})\sim t\M_{\on{inv}}$$
which is the equivalence of the two Cartan models. Conversely, 
suppose $\N$ is a $\g$-differential $W\g$-module. We have
$(S\g^*\otimes\N)_{\on{inv}}\sim \N_{\on{basic}}$ by Cartan's theorem
\ref{th:cargu}.
Therefore, the quasi-isomorphism from Theorem \ref{th:cartan}
$$ \N_{\on{basic}}\sim (S\g^*\otimes\N)_{\on{inv}}\sim t\N_{\on{inv}}$$
yields a quasi-isomorphism of differential $\wedge\P$-modules,
as in Theorem \ref{th:chevalley},
$$ h\N_{\on{basic}}\sim ht\N_{\on{inv}}\sim \N_{\on{inv}}.$$

\end{appendix}

\bibliographystyle{amsplain}   %

\def\polhk#1{\setbox0=\hbox{#1}{\ooalign{\hidewidth
  \lower1.5ex\hbox{`}\hidewidth\crcr\unhbox0}}} \def\cprime{$'$}
  \def\cprime{$'$} \def\polhk#1{\setbox0=\hbox{#1}{\ooalign{\hidewidth
  \lower1.5ex\hbox{`}\hidewidth\crcr\unhbox0}}} \def\cprime{$'$}
\providecommand{\bysame}{\leavevmode\hbox to3em{\hrulefill}\thinspace}
\providecommand{\MR}{\relax\ifhmode\unskip\space\fi MR }
\providecommand{\MRhref}[2]{%
  \href{http://www.ams.org/mathscinet-getitem?mr=#1}{#2}
}
\providecommand{\href}[2]{#2}

\end{document}